# Reiteration formulae for the real interpolation method including 𝓛 or 𝓡 limiting spaces.

**Leo R. Ya. Doktorski,**[1]

[1] Department Object Recognition, Fraunhofer Institute of Optronics, System Technologies and Image Exploitation IOSB, Gutleuthausstr. 1, 76275 Ettlingen, Germany.
Correspondence should be addressed to Leo R. Ya. Doktorski;
leo.doktorski@iosb.fraunhofer.de or doktorskileo@gmail.com.

### Abstract

We consider a real interpolation method defined by means of slowly varying functions. We present some reiteration formulae including so called 𝓛 or 𝓡 limiting interpolation spaces. These spaces arise naturally in reiteration formulae for the limiting cases $\theta = 0$ or $\theta = 1$. Applications to grand and small Lorentz spaces are given.

### 1. Introduction

Let $\bar{A} := (A_0, A_1)$ be a compatible couple of (quasi-) Banach spaces such that $A_0 \cap A_1 \neq \{0\}$. The Peetre's $K$-functional on $A_0 + A_1$ is given by

$$K(t,f;\bar{A}) \equiv K(t,f) := \inf\left(\|f_0\|_{A_0} + t\|f_1\|_{A_1}\right) \ (f = f_0 + f_1, f_i \in A_i, i=0,1, t>0). \quad (1)$$

The classical (Lions-Peetre) scale of interpolation spaces $\bar{A}_{\theta,q}$ is defined via the (quasi-)norms

$$\|f\|_{\theta,q} := \left\| u^{-\theta-1/q} K(u,f) \right\|_{q,(0,\infty)}, \quad (2)$$

where $\|*\|_{q,(a,b)}$ ($0<q\leq\infty$, $-\infty\leq a<b\leq\infty$) is the usual (quasi-)norm in the Lebesgue space $L_q$ on the interval $(a,b)$. For further information about basic properties of the $K$-functional and the real interpolation method we refer to e.g. [1].

Above definition requires $0<\theta<1$ or $\theta \in \{0,1\}$ for $q=\infty$, but limiting real interpolation spaces play an important role in certain questions in analysis. See for example [2 - 13] and references therein. One of the possibilities to consider the limiting cases $\theta \in \{0,1\}$ for $q < \infty$ is to involve an additional factor $a(u)$ in the formula (2):

$$\|f\|_{\theta,q;a} := \left\| u^{-\theta-1/q} a(u) K(u,f) \right\|_{q,(0,\infty)}.$$

For this factor various functions have been considered: logarithms [2 - 4], broken-logarithmic functions [5], or more generally a slowly varying functions [6 - 10] (see Definition 5 below). In this way one gets a scale $\bar{A}_{\theta,q;a}$, where the limiting cases $\theta \in \{0,1\}$ for $q < \infty$ have sense. An important property of this scale is that the reiteration spaces $\left(\bar{A}_{\theta_0,q_0,b_0}, \bar{A}_{\theta_1,q_1,b_1}\right)_{\theta,r,a}$ with $0 < \theta_0 < \theta_1 < 1$ and $0 < \theta < 1$ belong to the same scale. However, for limiting cases $\theta \in \{0,1\}$ the reiteration spaces normally do not belong to this scale and new interpolation functors are needed to describe them. Thus, the limiting reiteration formulae lead in a natural way to new interpolation spaces. Following [7], we call them 𝓡 und 𝓛 spaces (see Definition 6 below).

The principal aim of this paper is to establish reiteration formulae for couples where one or both operands are 𝓛 or 𝓡 spaces. This was motivated by articles [11 - 13], where it is shown that the so-called grand and small Lorentz spaces can be described in terms of the $\mathcal{R}$ and $\mathcal{L}$ spaces. Other motivation was the paper [9], which gives a good technical tool to prove Holmstedt's formulae for the $K$-functional, which are needed to establish the reiteration formulae.

This paper is organized as follows. Section 3 contains some notations, definitions and technical results. In Section 4, we collect necessary definitions and statements dealing with the real interpolation method involving slowly varying functions. In Sections 5 and 6, we establish reiteration formulae for the couples, where one of the operands is either $\mathcal{L}$ or $\mathcal{R}$ space. Reiteration formulae for the couples, where both operands are the $\mathcal{L}$ or $\mathcal{R}$ space are established in Section 7. In Section 8 is discussed how to obtain from our general theorems the ones for an ordered couple of (quasi-) Banach spaces. Finally, in Section 9, we present some interpolation results for the grand and small Lebesgue spaces as applications of our general reiteration theorems.

## 2. Preliminaries

We write $X \subset Y$ for two quasi-normed spaces $X$ and $Y$ to indicate that $X$ is continuously embedded in $Y$. The notation $X \cong Y$ means that $X \subset Y$ and $Y \subset X$. If $f$ and $g$ are positive functions, we write $f \prec g$ if $f \leq Cg$, where the constant $C$ is independent on all significant quantities. We say that two functions $f$ and $g$ are equivalent ($f \approx g$) if $f \prec g$ and $g \prec f$. We write $f \uparrow$ ($f \downarrow$) if the positive function $f$ is non-decreasing (non-increasing).

### 2.1. Slowly varying functions

We will need definition and some basic properties of slowly varying functions. For further information about slowly varying functions we refer to e.g. [6 - 8].

**Definition 1.** We say that a positive Lebesgue measurable function $b$ is slowly varying on $(0,\infty)$ (or on $(0,1)$), notation $b \in SV(0,\infty) \equiv SV$ (or $b \in SV(0,1)$), if, for each $\varepsilon > 0$, the function $t^{\varepsilon}b(t)$ is equivalent to an increasing function while the function $t^{-\varepsilon}b(t)$ is equivalent to a decreasing function.

**Lemma 2.** Let $b, b_1, b_2 \in SV$, $r>0$, $\alpha>0$, $0<q\leq\infty$, and $t \in (0, \infty)$.
  (i)  Then $b^r \in SV$, $b\left(\frac{1}{t}\right) \in SV$, $b(t^r b_1(t)) \in SV$, and $b_1 b_2 \in SV$.
  (ii) If $f \approx g$, then $b \circ f \approx b \circ g$.
  (iii) $\left\|u^{\alpha-1/q}b(u)\right\|_{q,(0,t)} \approx t^{\alpha}b(t)$ and $\left\|u^{-\alpha-1/q}b(u)\right\|_{q,(t,\infty)} \approx t^{-\alpha}b(t)$.
  (iv) The functions $\left\|u^{-1/q}b(u)\right\|_{q,(0,t)}$ and $\left\|u^{-1/q}b(u)\right\|_{q,(t,\infty)}$ (if exist) belong to $SV$ and $b(t) \prec \left\|u^{-1/q}b(u)\right\|_{q,(0,t)}$ and $b(t) \prec \left\|u^{-1/q}b(u)\right\|_{q,(t,\infty)}$.

**Remark 3.** Let $a, b \in SV$, $\lambda > 0$, and $\rho(t) = t^{\lambda}a(t)$.
  (i) We may assume without loss of generality that the function $a$ is continuous.
  (ii) Due to [10, Lemma 3.1], there exists a strongly increasing, differentiable function $\sigma(t) \approx \rho(t)$, such that $\sigma(+0) = 0$, $\sigma(\infty) = \infty$, and
$$\sigma'(t) \approx \sigma(t)\frac{1}{t}. \tag{3}$$
Obviously, the function $\sigma$ has inverse $\sigma^{(-1)}$ and it holds
$$\sigma^{(-1)'}(t) \approx \sigma^{(-1)}(t)\frac{1}{t}. \tag{4}$$
  (iii) It is not difficult to show that for any compatible couple $\bar{A}$, it holds

$$K(\rho(t), f; \bar{A}) \approx K(\sigma(t), f; \bar{A}) \tag{5}$$

for all $f \in A_0 + A_1$ and $t \in (0, \infty)$.

(iv) Using [6, Proposition 2.2 (ii)], it can be shown that $b \circ \sigma \in SV$ and $b \circ \sigma^{(-1)} \in SV$.

(v) Let $0 \leq \theta \leq 1$, $0 < r \leq \infty$. Consider the expression

$$\left\| \rho(t)^{-\theta} t^{-1/r} b(\rho(t)) \bar{K}(\rho(t), f) \right\|_{r,(0,\infty)}.$$

(5) implies that for all $f \in A_0 + A_1$

$$\left\| \rho(t)^{-\theta} t^{-1/r} a(\rho(t)) \bar{K}(\rho(t), f) \right\|_{r,(0,\infty)} \approx \left\| \sigma(t)^{-\theta} t^{-1/r} a(\sigma(t)) \bar{K}(\sigma(t), f) \right\|_{r,(0,\infty)}.$$

Using change of variables $u = \sigma(t)$ and (3), we get

$$\left\| \sigma(t)^{-\theta} t^{-1/r} a(\sigma(t)) \bar{K}(\sigma(t), f) \right\|_{r,(0,\infty)} \approx \left\| u^{-\theta-1/r} a(u) \bar{K}(u, f) \right\|_{r,(0,\infty)}.$$

Thus,

$$\left\| \rho(t)^{-\theta} t^{-1/r} b(\rho(t)) \bar{K}(\rho(t), f) \right\|_{r,(0,\infty)} \approx \left\| u^{-\theta-1/r} a(u) \bar{K}(u, f) \right\|_{r,(0,\infty)}.$$

We will use the last property when making the change of variables.

### 2.2. Hardy-type inequalities

We need the following Hardy-type inequalities [6, Lemma 2.7].

**Lemma 4.** Let $1 \leq P \leq \infty$ and $b \in SV$.

(i) The inequality

$$\left\| t^{\nu-1/P} b(t) \int_0^t g(u) du \right\|_{P,(0,\infty)} \prec \left\| t^{\nu+1-1/P} b(t) g(t) \right\|_{P,(0,\infty)}$$

holds for all (Lebesgue-) measurable nonnegative functions $g$ on $(0,\infty)$ if and only if $\nu < 0$.

(ii) The inequality

$$\left\| t^{\nu-1/P} b(t) \int_t^\infty g(u) du \right\|_{P,(0,\infty)} \prec \left\| t^{\nu+1-1/P} b(t) g(t) \right\|_{P,(0,\infty)}$$

holds for all (Lebesgue-) measurable nonnegative functions $g$ on $(0,\infty)$ if and only if $\nu > 0$.

## 3. Real interpolation method involving slowly varying functions: a short overview

Here we collect necessary definitions and statements dealing with the real interpolation method involving slowly varying functions. We consider a compatible couple of (quasi-) Banach spaces $\bar{A} = (A_0, A_1)$ such that $A_0 \cap A_1 \neq \{0\}$ and use the Peetre's $K$-functional on $A_0 + A_1$, given by (1).

### 3.1. Standard interpolation spaces

**Definition 5.** ([6]). Let $0 \leq \theta \leq 1$, $0 < q \leq \infty$ and $b \in SV$. We write

$$\bar{A}_{\theta,q;b} \equiv (A_0, A_1)_{\theta,q;b} := \left\{ f \in A_0 + A_1 : \|f\|_{\theta,q;b} = \left\| u^{-\theta-1/q} b(u) K(u, f) \right\|_{q,(0,\infty)} < \infty \right\}.$$

It is known [6], that the spaces $\bar{A}_{\theta,q;b}$ are intermediate spaces between $A_0$ and $A_1$ if and only if one of the following conditions is satisfied:

(i) $0 < \theta < 1$,

(ii) $\theta = 0$ and $\left\| u^{-1/q} b(u) \right\|_{q,(1,\infty)} < \infty$,

(iii) $\theta = 1$ and $\left\|u^{-1/q} b(u)\right\|_{q,(0,1)} < \infty$.

We refer to the spaces $\bar{A}_{\theta,q;b}$ as standard interpolation spaces.

### 3.2. $\mathcal{R}$ and $\mathcal{L}$ limiting interpolation spaces

**Definition 6.** ([6]). Let $0 < r, q \leq \infty$, $0 < \sigma < 1$, $a, b \in SV$. We write

$$\bar{A}^{\mathcal{L}}_{\sigma,r,b,q,a} \equiv (A_0, A_1)^{\mathcal{L}}_{\sigma,r,b,q,a} :=$$

$$\left\{ f \in A_0 + A_1 : \|f\|_{\mathcal{L};\sigma,r,b,q,a} := \left\| t^{-1/r} \frac{b(t)}{a(t)} \left\| u^{-\sigma-1/q} a(u) K(u,f) \right\|_{q,(0,t)} \right\|_{r,(0,\infty)} < \infty \right\}.$$

Similarly,

$$\bar{A}^{\mathcal{R}}_{\sigma,r,b,q,a} \equiv (A_0, A_1)^{\mathcal{R}}_{\sigma,r,b,q,a} :=$$

$$\left\{ f \in A_0 + A_1 : \|f\|_{\mathcal{R};\sigma,r,b,q,a} := \left\| t^{-1/r} \frac{b(t)}{a(t)} \left\| u^{-\sigma-1/q} a(u) K(u,f) \right\|_{q,(t,\infty)} \right\|_{r,(0,\infty)} < \infty \right\}.$$

It is known [6], that the spaces $\bar{A}^{\mathcal{L}}_{\sigma,r,b,q,a}$ ($\bar{A}^{\mathcal{R}}_{\sigma,r,b,q,a}$) are intermediate spaces between $A_0$ and $A_1$ if and only if $\left\| s^{-1/r} \frac{b(s)}{a(s)} \right\|_{r,(1,\infty)} < \infty$ ($\left\| s^{-1/r} \frac{b(s)}{a(s)} \right\|_{r,(0,1)} < \infty$, respectively).

Similar definitions can be found for example in [2, 3, 5, 7, 8, 10]. In those papers the reader can find further properties of these spaces. Following [7], we refer to the spaces $\bar{A}^{\mathcal{L}}_*$ and $\bar{A}^{\mathcal{R}}_*$ as $\mathcal{L}$ and $\mathcal{R}$ spaces (or $\mathcal{L}$ and $\mathcal{R}$ limiting interpolation spaces).

### 3.3. Some known formulae and reiteration theorems

Under suitable conditions, the following formulae hold [1, Chap. 5, Proposition 1.2], [6]:

$$K(t, f; A_0, A_1) = t K(t^{-1}, f; A_1, A_0), t > 0,$$

$$(A_0, A_1)_{\theta,q,b} = (A_1, A_0)_{1-\theta,q,\tilde{b}}, \tag{6}$$

$$(A_0, A_1)^{\mathcal{L}}_{\sigma,r,b,q,a} = (A_1, A_0)^{\mathcal{R}}_{1-\sigma,r,\tilde{b},q,\tilde{a}}, \qquad \text{where } \tilde{b}(t) = b(t^{-1}). \tag{7}$$

The following reiteration theorems in different modifications can be found e.g. in [6, 8, 10, 13]. We formulate them in the form we will use them.

**Theorem 7.** ([6, (3.4)]). Let $0 < \theta < 1$, $0 < q_0, q_1, r \leq \infty$, $a, b_0, b_1 \in SV$, $0 < \theta_0 < \theta_1 < 1$. Put $\rho(t) = t^{\theta_1 - \theta_0} \frac{b_0(t)}{b_1(t)}$. Then

$$\left( \bar{A}_{\theta_0,q_0,b_0}, \bar{A}_{\theta_1,q_1,b_1} \right)_{\theta,r,a} \cong \bar{A}_{\eta,r,a^\#},$$

where $\eta = (1-\theta)\theta_0 + \theta\theta_1$ and $a^\# = b_0^{1-\theta} b_1^\theta a \circ \rho$.

**Theorem 8.** (Cf. [8, Theorems 5.10 and 5.12], [13, Lemma 5.1.], [10].) Let $0 < q, r \leq \infty$ and $a, b \in SV$.

(i) If $\left\| t^{-1/q} b(t) \right\|_{q,(1,\infty)} < \infty$ and $0 < \theta \leq 1$, then

$$\left(\bar{A}_{0,q,b}, A_1\right)_{\theta,r,a} \cong \bar{A}_{\theta,r,a^{\#}},$$

where $a^{\#}(t) = \left(\|s^{-1/q}b(s)\|_{q,(t,\infty)}\right)^{1-\theta} a\left(t\|s^{-1/q}b(s)\|_{q,(t,\infty)}\right)$.

(ii) If $\|t^{-1/q}b(t)\|_{q,(0,1)} < \infty$ and $0 \leq \theta < 1$, then
$$\left(A_0, \bar{A}_{1,q,b}\right)_{\theta,r,a} \cong \bar{A}_{\theta,r,a^{\#}},$$

where $a^{\#}(t) = \left(\|s^{-1/q}b(s)\|_{q,(0,t)}\right)^{\theta} a\left(t/\|s^{-1/q}b(s)\|_{q,(0,t)}\right)$.

**Theorem 9.** Let $0<\theta_0<1$, $0<q_0, q_1, r \leq \infty$, and $a, b_0, b_1 \in SV$.

(i) (Cf. [7, Theorem 5.4] and [6, (4.13)].) If $0<\theta \leq 1$, then
$$\left(\bar{A}_{\theta_0,q_0,b_0}, A_1\right)_{\theta,r,a} \cong \bar{A}_{\eta,r,a^{\#}},$$

where $\eta = (1-\theta)\theta_0 + \theta$, $a^{\#}(t) = b_0(t)^{1-\theta}a(\rho(t))$, and $\rho(t) = t^{1-\theta_0}b_0(t)$

(ii) (Cf. [8, Theorem 5.2].) If $0<\theta<1$ and $\|s^{-1/q_1}b_1(s)\|_{q_1,(0,1)} < \infty$, then
$$\left(\bar{A}_{\theta_0,q_0,b_0}, \bar{A}_{1,q_1,b_1}\right)_{\theta,r,a} \cong \bar{A}_{\eta,r,a^{\#}},$$

where $\eta = (1-\theta)\theta_0 + \theta$, $a^{\#}(t) = b_0(t)^{1-\theta}\left(\|s^{-1/q_1}b_1(s)\|_{q_1,(0,t)}\right)^{\theta} a(\rho(t))$, and $\rho(t) = t^{1-\theta_0}\frac{b_0(t)}{\|s^{-1/q_1}b_1(s)\|_{q_1,(0,t)}}$.

**Theorem 10.** (Cf. [6, (3.21) and Theorem 3.2] and [8, (5.4)].) Suppose that $0 < \theta_0 < \theta_1 < 1$, $0<r_0, q_0, q_1 \leq \infty$, $a_0, a_1, b_0 \in SV$, and $\left\|s^{-1/r_0}\frac{b_0(s)}{a_0(s)}\right\|_{r_0,(1,\infty)} < \infty$.

(i) Let $\sigma$ be a strongly increasing, differentiable function such that $\sigma(t) \approx t^{1-\theta_0}a_0(t)$ and $a^{\#} = \frac{b_0 \circ \sigma^{(-1)}}{a_0 \circ \sigma^{(-1)}}$. Then
$$\left(\bar{A}_{\theta_0,q_0,a_0}, A_1\right)_{0,r_0,a^{\#}} \cong \bar{A}^{\mathcal{L}}_{\theta_0,r_0,b_0,q_0,a_0}.$$

(ii) Let $\|s^{-1/q_1}a_1(s)\|_{q_1,(0,1)} < \infty$, $\sigma$ be a strongly increasing, differentiable function such that $\sigma(t) \approx t^{1-\theta_0}\frac{a_0(t)}{\|s^{-1/q_1}a_1(s)\|_{q_1,(0,t)}}$ and $a^{\#} = \frac{b_0 \circ \sigma^{(-1)}}{a_0 \circ \sigma^{(-1)}}$. Then
$$\left(\bar{A}_{\theta_0,q_0,a_0}, \bar{A}_{1,q_1,a_1}\right)_{0,r_0,a^{\#}} \cong \bar{A}^{\mathcal{L}}_{\theta_0,r_0,b_0,q_0,a_0}.$$

(iii) Let $\sigma$ be a strongly increasing, differentiable function such that $\sigma(t) \approx t^{\theta_1-\theta_0}\frac{a_0(t)}{a_1(t)}$ and $a^{\#} = \frac{b_0 \circ \sigma^{(-1)}}{a_0 \circ \sigma^{(-1)}}$. Then
$$\left(\bar{A}_{\theta_0,q_0,a_0}, \bar{A}_{\theta_1,q_1,a_1}\right)_{0,r_0,a^{\#}} \cong \bar{A}^{\mathcal{L}}_{\theta_0,r_0,b_0,q_0,a_0}.$$

*Proof.* We prove only the first statement. Let $\rho(t) = t^{1-\theta_0}a_0(t)$ and $\sigma \approx \rho$. Note that such function $\sigma$ exists (see Remark 3). By [6, (3.21)] we have to show that
$$(a_0 a^{\#} \circ \rho) \approx b_0.$$

Due to Lemma 2, $a^{\#} \circ \sigma \in SV$ and $a^{\#} \circ \rho \approx a^{\#} \circ \sigma$. Thus
$$(a_0 a^{\#} \circ \rho)(t) \approx (a_0 a^{\#} \circ \sigma)(t) = a_0(t)\frac{b_0\left(\sigma^{(-1)}(\sigma(t))\right)}{a_0\left(\sigma^{(-1)}(\sigma(t))\right)} = b_0(t).$$

The statements (ii) and (iii) can be proved analogously. □

## 4. Interpolation between the $\mathcal{L}$ spaces and the standard interpolation spaces

At first, we establish reiteration formulae for the couples in which the $\mathcal{L}$ space is the first operand.

**Theorem 11.** Let $0 < \theta_0 < \theta_1 < 1$, $0 < q_0, r_0, r_1, r \leq \infty$, $a, a_0, b_0, b_1 \in SV$, and $\left\|s^{-1/r_0}\frac{b_0(s)}{a_0(s)}\right\|_{r_0,(1,\infty)} < \infty$. Put $c_0(t) = a_0(t)\left\|s^{-1/r_0}\frac{b_0(s)}{a_0(s)}\right\|_{r_0,(t,\infty)}$ and $\rho(t) = t^{\theta_1-\theta_0}\frac{c_0(t)}{b_1(t)}$.

(i) If $0 < \theta < 1$ then
$$\left(\bar{A}^{\mathcal{L}}_{\theta_0,r_0,b_0,q_0,a_0}, \bar{A}_{\theta_1,r_1,b_1}\right)_{\theta,r,a} \cong \bar{A}_{\eta,r,a^\#},$$
where $\eta = (1-\theta)\theta_0 + \theta\theta_1$ and $a^\#(t) = c_0(t)^{1-\theta}b_1(t)^\theta a(\rho(t))$.

(ii) If additionally $\left\|s^{-1/r}a(s)\right\|_{r,(0,1)} < \infty$, then
$$\left(\bar{A}^{\mathcal{L}}_{\theta_0,r_0,b_0,q_0,a_0}, \bar{A}_{\theta_1,r_1,b_1}\right)_{1,r,a} \cong \vec{A}^{\mathcal{R}}_{\theta_1,r,a^\#,r_1,b_1},$$
where $a^\#(t) = b_1(t)a(\rho(t))$.

*Proof.* Denote $X_0 = \bar{A}_{\theta_0,q_0,a_0}$, $X_1 = \bar{A}_{\theta_1,r_1,b_1}$, and $Z = \left(\bar{A}^{\mathcal{L}}_{\theta_0,r_0,b_0,q_0,a_0}, \bar{A}_{\theta_1,r_1,b_1}\right)_{\theta,r,a}$. Let $\sigma$ be a strongly increasing, differentiable function such that $\sigma(t) \approx t^{\theta_1-\theta_0}\frac{a_0(t)}{b_1(t)}$ (see Remark 3 (ii)) and $c = \frac{b_0 \circ \sigma^{(-1)}}{a_0 \circ \sigma^{(-1)}}$. By Theorem 10 (iii) we obtain
$$\vec{A}^{\mathcal{L}}_{\theta_0,r_0,b_0,q_0,a_0} \cong \left(\bar{A}_{\theta_0,q_0,a_0}, \bar{A}_{\theta_1,r_1,b_1}\right)_{0,r_0,c} = \bar{X}_{0,r_0,c},$$
Hence, $Z \cong \left(\bar{X}_{0,r_0,c}, X_1\right)_{\theta,r,a}$. Using Theorem 8 (i), we conclude that for $0 < \theta \leq 1$, it holds
$$Z \cong \bar{X}_{\theta,r,d} = \left(\bar{A}_{\theta_0,q_0,a_0}, \bar{A}_{\theta_1,r_1,b_1}\right)_{\theta,r,d}, \tag{8}$$
where $d(t) = \left(\left\|s^{-1/r_0}c(s)\right\|_{r_0,(t,\infty)}\right)^{1-\theta} a\left(t\left\|s^{-1/r_0}c(s)\right\|_{r_0,(t,\infty)}\right)$.

By change of variables $u = \sigma^{(-1)}(s)$ (see Remark 3) we arrive at
$$\left\|s^{-1/r_0}c(s)\right\|_{r_0,(\sigma(t),\infty)} \approx \left\|u^{-1/r_0}\frac{b_0(u)}{a_0(u)}\right\|_{r_0,(t,\infty)}.$$
Thus, because $\sigma(t)\left\|u^{-1/r_0}\frac{b_0(u)}{a_0(u)}\right\|_{r_0,(t,\infty)} \approx \rho(t)$,
$$d(\sigma(t)) \approx \left(\left\|u^{-1/r_0}\frac{b_0(u)}{a_0(u)}\right\|_{r_0,(t,\infty)}\right)^{1-\theta} a(\rho(t)). \tag{9}$$

(i) Let $0 < \theta < 1$. By (8) and Theorem 7 it follows $Z \cong \bar{A}_{\eta,r,a^\#}$, where $a^\#(t) = a_0(t)^{1-\theta}b_1(t)^\theta d(\sigma(t))$. Hence by (9),
$$a^\#(t) \approx \left(a_0(t)\left\|u^{-1/r_0}\frac{b_0(u)}{a_0(u)}\right\|_{r_0,(t,\infty)}\right)^{1-\theta} b_1(t)^\theta a(\rho(t)) = c_0(t)^{1-\theta}b_1(t)^\theta a(\rho(t)).$$

(ii) Let $\theta = 1$. From (8) we have
$$\left(\bar{A}^{\mathcal{L}}_{\theta_0,r_0,b_0,q_0,a_0}, \bar{A}_{\theta_1,r_1,b_1}\right)_{1,r,a} \cong \left(\bar{A}_{\theta_0,q_0,a_0}, \bar{A}_{\theta_1,r_1,b_1}\right)_{1,r,d},$$
where $d(t) = a\left(t\left\|s^{-1/r_0}c(s)\right\|_{r_0,(t,\infty)}\right)$. By [6, (4.7)] we conclude that
$$\left(\bar{A}_{\theta_0,q_0,a_0}, \bar{A}_{\theta_1,r_1,b_1}\right)_{1,r,d} \cong \vec{A}^{\mathcal{R}}_{\theta_1,r,a^\#,r_1,b_1},$$
where by (9),
$$a^\#(t) = b_1(t)d(\sigma(t)) \approx b_1(t)a(\rho(t)).$$

This completes the proof. □

**Theorem 12.** Let $0<\theta,\theta_0<1$, $0 < q_0, r_0, r_1, r \leq \infty$, $a, a_0, b_0, b_1 \in SV$, $\left\|s^{-1/r_0}\frac{b_0(s)}{a_0(s)}\right\|_{r_0,(1,\infty)} < \infty$, and $\left\|s^{-1/r_1}b_1(s)\right\|_{r_1,(0,1)} < \infty$. Put $c_0(t) = a_0(t)\left\|s^{-1/r_0}\frac{b_0(s)}{a_0(s)}\right\|_{r_0,(t,\infty)}$, $c_1(t) = \left\|s^{-1/r_1}b_1(s)\right\|_{r_1,(0,t)}$, and $\rho(t) = t^{1-\theta_0}\frac{c_0(t)}{c_1(t)}$. Then
$$\left(\vec{A}^{\mathcal{L}}_{\theta_0,r_0,b_0,q_0,a_0}, \bar{A}_{1,r_1,b_1}\right)_{\theta,r,a} \cong \bar{A}_{\eta,r,a^{\#}},$$
where $\eta = (1-\theta)\theta_0 + \theta$ and $a^{\#}(t) = c_0(t)^{1-\theta}c_1(t)^{\theta}a(\rho(t))$.

*Proof.* Denote $X_0 = \bar{A}_{\theta_0,q_0,a_0}$, $X_1 = \bar{A}_{1,r_1,b_1}$, and $Z = \left(\vec{A}^{\mathcal{L}}_{\theta_0,r_0,b_0,q_0,a_0}, \bar{A}_{1,r_1,b_1}\right)_{\theta,r,a}$. Let $\sigma$ be a strongly increasing, differentiable function such that $\sigma(t) \approx t^{1-\theta_0}\frac{a_0(t)}{\left\|s^{-1/r_1}b_1(s)\right\|_{r_1,(0,t)}}$ and $c = \frac{b_0 \circ \sigma^{(-1)}}{a_0 \circ \sigma^{(-1)}}$. By Theorem 10 (ii) we have
$$\vec{A}^{\mathcal{L}}_{\theta_0,r_0,b_0,q_0,a_0} \cong \left(\bar{A}_{\theta_0,q_0,a_0}, \bar{A}_{1,r_1,b_1}\right)_{0,r_0,c} = \bar{X}_{0,r_0,c}.$$
Hence, $Z \cong \left(\bar{X}_{0,r_0,c}, X_1\right)_{\theta,r,a}$. Using Theorem 8 (i), we conclude that
$$Z \cong \bar{X}_{\theta,r,d} = \left(\bar{A}_{\theta_0,q_0,a_0}, \bar{A}_{1,r_1,b_1}\right)_{\theta,r,d},$$
where $d(t) = \left(\left\|s^{-1/r_0}c(s)\right\|_{r_0,(t,\infty)}\right)^{1-\theta}a\left(t\left\|s^{-1/r_0}c(s)\right\|_{r_0,(t,\infty)}\right)$. By Theorem 9 (ii) we derive that
$$Z \cong \left(\bar{A}_{\theta_0,q_0,a_0}, \bar{A}_{1,r_1,b_1}\right)_{\theta,r,d} \cong \bar{A}_{\eta,r,a^{\#}},$$
where by (9)
$$a^{\#}(t) = a_0(t)^{1-\theta}c_1(t)^{\theta}d(\sigma(t)) =$$
$$= \left(a_0(t)\left\|s^{-1/r_0}c(s)\right\|_{r_0,(\sigma(t),\infty)}\right)^{1-\theta}c_1(t)^{\theta}a\left(\sigma(t)\left\|s^{-1/r_0}c(s)\right\|_{r_0,(\sigma(t),\infty)}\right) \approx$$
$$\approx c_0(t)^{1-\theta}c_1(t)^{\theta}a(\rho(t)).$$
This completes the proof. □

The following theorem can be proved analogously.

**Theorem 13.** Let $0<\theta_0<1$, $0 < \theta \leq 1$, $0<q_0, r_0, r\leq\infty$, $a, a_0, b_0 \in SV$, and $\left\|s^{-1/r_0}\frac{b_0(s)}{a_0(s)}\right\|_{r_0,(1,\infty)} < \infty$. If $\theta = 1$, we additionally assume that $\left\|s^{-1/r}a(s)\right\|_{r,(0,1)} < \infty$. Put $c_0(t) = a_0(t)\left\|s^{-1/r_0}\frac{b_0(u)}{a_0(u)}\right\|_{r_0,(t,\infty)}$ and $(t) = t^{1-\theta_0}c_0(t)$. Then
$$\left(\vec{A}^{\mathcal{L}}_{\theta_0,r_0,b_0,q_0,a_0}, A_1\right)_{\theta,r,a} \cong \bar{A}_{\eta,r,a^{\#}},$$
where $\eta = (1-\theta)\theta_0 + \theta$ and $a^{\#}(t) = c_0(t)^{1-\theta}a(\rho(t))$.

Now we establish reiteration formulae for the couples in which the $\mathcal{L}$ space is the second operand. Here we need corresponding Holmstedt-type formulae (Theorem 14 and Theorem 17), which will be proved, based on the theorems from [9]. Because of this we will adopt the notation from that paper.

**Theorem 14.** Let $0 < \theta_1 < 1$, $0 < q_1, r_1 \leq \infty$, $a_1, b_1 \in SV$, and $\left\|s^{-1/r_1}\frac{b_1(s)}{a_1(s)}\right\|_{r_1,(1,\infty)} < \infty$.
Put $c_1(t) = a_1(t)\left\|s^{-1/r_1}\frac{b_1(s)}{a_1(s)}\right\|_{r_1,(t,\infty)}$ and $\rho(t) = t^{\theta_1}\frac{1}{c_1(t)}$. Then for all $f \in A_0 +$
$\bar{A}^{\mathcal{L}}_{\theta_1,r_1,b_1,q_1,a_1}$ and $t > 0$,
$$K\left(\rho(t), f; A_0, \bar{A}^{\mathcal{L}}_{\theta_1,r_1,b_1,q_1,a_1}\right) \approx$$
$$\approx K(t,f) + \rho(t)\left\|s^{-1/r_1}\frac{b_1(s)}{a_1(s)}\left\|u^{-\theta_1-1/q_1}a_1(u)K(u,f)\right\|_{q_1,(t,s)}\right\|_{r_1,(t,\infty)}.$$

*Proof.* We will often be using the condition $0 < \theta_1 < 1$ and the basic properties of slowly varying functions formulated in Lemma 2. Let $\Phi_1$ be the function space corresponding to $\bar{A}^{\mathcal{L}}_{\theta_1,r_1,b_1,q_1,a_1}$:
$$\|F(*)\|_{\Phi_1} = \left\|s^{-1/r_1}\frac{b_1(s)}{a_1(s)}\left\|u^{-\theta_1-1/q_1}a_1(u)F(u)\right\|_{q_1,(0,s)}\right\|_{r_1,(0,\infty)}.$$
Consider the functions $J(t,f)$, $g_1(t)$, and $h_1(t)$.
$$J(t,f) := \left\|\chi_{(t,\infty)}(*)K(*,f)\right\|_{\Phi_1} = \left\|s^{-1/r_1}\frac{b_1(s)}{a_1(s)}\left\|u^{-\theta_1-1/q_1}a_1(u)K(u,f)\right\|_{q_1,(t,s)}\right\|_{r_1,(t,\infty)}.$$
$$g_1(t) := t\left\|\chi_{(t,\infty)}(*)\right\|_{\Phi_1} = t\left\|s^{-1/r_1}\frac{b_1(s)}{a_1(s)}\left\|u^{-\theta_1-1/q_1}a_1(u)\right\|_{q_1,(t,s)}\right\|_{r_1,(t,\infty)} \leq$$
$$\leq t\left\|s^{-1/r_1}\frac{b_1(s)}{a_1(s)}\left\|u^{-\theta_1-1/q_1}a_1(u)\right\|_{q_1,(t,\infty)}\right\|_{r_1,(t,\infty)} \approx t^{1-\theta_1}c_1(t).$$
$$h_1(t) := \left\|*\chi_{(0,t)}(*)\right\|_{\Phi_1} = \left\|s^{-1/r_1}\frac{b_1(s)}{a_1(s)}\left\|u^{1-\theta_1-1/q_1}a_1(u)\right\|_{q_1,(0,\min(t,s))}\right\|_{r_1,(0,\infty)} \approx$$
$$\approx \left\|s^{-1/r_1}\frac{b_1(s)}{a_1(s)}(\min(t,s))^{1-\theta_1}a_1(\min(t,s))\right\|_{r_1,(0,\infty)} \approx$$
$$\approx \left\|s^{1-\theta_1-1/r_1}b_1(s)\right\|_{r_1,(0,t)} + t^{1-\theta_1}a_1(t)\left\|s^{-1/r_1}\frac{b_1(s)}{a_1(s)}\right\|_{r_1,(t,\infty)} \approx$$
$$\approx t^{1-\theta_1}b_1(t) + t^{1-\theta_1}c_1(t).$$
Because $b_1(t) = a_1(t)\frac{b_1(t)}{a_1(t)} \prec a_1(t)\left\|s^{-1/r_1}\frac{b_1(s)}{a_1(s)}\right\|_{r_1,(t,\infty)} = c_1(t)$, we get
$$g_1(t) \prec h_1(t) \approx t^{1-\theta_1}c_1(t).$$
Thus, we arrive at
$$g_1(t) + h_1(t) \approx t^{1-\theta_1}c_1(t) = \frac{t}{\rho(t)}.$$
Hence, by [9, Theorem 3], it yields
$$K\left(\rho(t), f; A_0, \bar{A}^{\mathcal{L}}_{\theta_1,r_1,b_1,q_1,a_1}\right) \approx$$
$$\approx K(t,f) + \rho(t)\left\|s^{-1/r_1}\frac{b_1(s)}{a_1(s)}\left\|u^{-\theta_1-1/q_1}a_1(u)K(u,f)\right\|_{q_1,(t,s)}\right\|_{r_1,(t,\infty)}.$$
This completes the proof. □

**Theorem 15.** Let $0 < \theta_1 < 1$, $0 \leq \theta < 1$, $0 < q_1, r, r_1 \leq \infty$, $a, a_1, b_1 \in SV$, and $\left\|s^{-1/r_1}\frac{b_1(s)}{a_1(s)}\right\|_{r_1,(1,\infty)} < \infty$. If $\theta = 0$, we additionally assume that $\left\|s^{-1/r}a(s)\right\|_{r,(1,\infty)} < \infty$.
Put $c_1(t) = a_1(t)\left\|s^{-1/r_1}\frac{b_1(s)}{a_1(s)}\right\|_{r_1,(t,\infty)}$ and $\rho(t) = t^{\theta_1}\frac{1}{c_1(t)}$. Then
$$\left(A_0, \bar{A}^{\mathcal{L}}_{\theta_1,r_1,b_1,q_1,a_1}\right)_{\theta,r,a} \cong \bar{A}_{\eta,r,a^{\#}},$$
where $\eta = \theta\theta_1$ and $a^{\#}(t) = c_1(t)^{\theta}a(\rho(t))$.

*Proof.* Denote $X_0 = A_0$, $X_1 = \bar{A}^{\mathcal{L}}_{\theta_1,r_1,b_1,q_1,a_1}$, $\bar{K}(t,f) = K(t,f;\bar{X})$, and $Z = \left(A_0, \bar{A}^{\mathcal{L}}_{\theta_1,r_1,b_1,q_1,a_1}\right)_{\theta,r,a}$. By Remark 3 (v) and Theorem 14 we can write

$$\|f\|_Z = \left\|t^{-\theta-1/r}a(t)\overline{K}(u,f)\right\|_{r,(0,\infty)} \approx$$
$$\approx \left\|\rho(t)^{-\theta}t^{-1/r}a(\rho(t))\overline{K}(\rho(t),f)\right\|_{r,(0,\infty)} =$$
$$= \left\|t^{-\eta-1/r}c_1(t)^\theta a(\rho(t))\overline{K}(\rho(t),f)\right\|_{r,(0,\infty)} \approx I_1 + I_2,$$

where
$$I_1 := \left\|t^{-\eta-1/r}c_1(t)^\theta a(\rho(t))K(t,f)\right\|_{r,(0,\infty)} = \|f\|_{\eta,r;a^\#}$$

and
$$I_2 :=$$
$$\left\|t^{\theta_1-\eta-1/r}c_1(t)^{\theta-1}a(\rho(t))\left\|s^{-1/r_1}\frac{b_1(s)}{a_1(s)}\left\|u^{-\theta_1-1/q_1}a_1(u)K(u,f)\right\|_{q_1,(t,s)}\right\|_{r_1,(t,\infty)}\right\|_{r,(0,\infty)}.$$

The proof of the theorem is completed, if we show that $I_2 \prec I_1$. We have
$$\left\|s^{-1/r_1}\frac{b_1(s)}{a_1(s)}\left\|u^{-\theta_1-1/q_1}a_1(u)K(u,f)\right\|_{q_1,(t,s)}\right\|_{r_1,(t,\infty)} \leq$$
$$\left\|s^{-1/r_1}\frac{b_1(s)}{a_1(s)}\right\|_{r_1,(t,\infty)}\left\|u^{-\theta_1-1/q_1}a_1(u)K(u,f)\right\|_{q_1,(t,\infty)}.$$

Denote $d(t) = c_1(t)^{\theta-1}a(\rho(t))\left\|s^{-1/r_1}\frac{b_1(s)}{a_1(s)}\right\|_{r_1,(t,\infty)}$. Then we get
$$I_2 \leq \left\|t^{\theta_1-\eta-1/r}d(t)\left\|u^{-\theta_1-1/q_1}a_1(u)K(u,f)\right\|_{q_1,(t,\infty)}\right\|_{r,(0,\infty)} := I_3.$$

It is enough to show that $I_3 \prec I_1$.

Case $q_1 = r$. Because $\theta_1 - \eta > 0$ and $d(u)a_1(u) = c_1(t)^\theta a(\rho(t))$, we have
$$I_3 = \left\|t^{\theta_1-\eta-1/r}d(t)\left\|u^{-\theta_1-1/r}a_1(u)K(u,f)\right\|_{r,(t,\infty)}\right\|_{r,(0,\infty)} =$$
$$= \left\|\left\|t^{\theta_1-\eta-1/r}d(t)\right\|_{r,(0,u)}u^{-\theta_1-1/r}a_1(u)K(u,f)\right\|_{r,(0,\infty)} \approx$$
$$\approx \left\|u^{-\eta-1/r}d(u)a_1(u)K(u,f)\right\|_{r,(0,\infty)} = I_1.$$

Case $q_1 < r$. In this case, using Lemma 4 (ii) (Hardy-type inequality), we also get $I_3 \prec I_1$. The details are left to the reader.

Case $r < q_1$. We need the following estimate:
$$\left\|u^{-\theta_1-1/q_1}a_1(u)K(u,f)\right\|_{q_1,(t,\infty)} \prec \left\|u^{-\theta_1-1/r}a_1(u)K(u,f)\right\|_{r,(t,\infty)}. \tag{10}$$

For $u>t$, because $K(u,f)\uparrow$, we have
$$u^{-\theta_1}a_1(u)K(u,f) \approx K(u,f)\left\|s^{-\theta_1-1/r}a_1(s)\right\|_{r,(u,\infty)} \leq$$
$$\leq \left\|s^{-\theta_1-1/r}a_1(s)K(s,f)\right\|_{r,(u,\infty)} \leq \left\|s^{-\theta_1-1/r}a_1(s)K(s,f)\right\|_{r,(t,\infty)}.$$

Thus, (10) is proved for $q_1=\infty$:
$$\sup_{t<u<\infty} u^{-\theta_1}a_1(u)K(u,f) \prec \left\|u^{-\theta_1-1/r}a_1(u)K(u,f)\right\|_{r,(t,\infty)}.$$

If $r<q_1<\infty$ using the last estimate, we get
$$\left\|u^{-\theta_1-1/q_1}a_1(u)K(u,f)\right\|_{q_1,(t,\infty)} = \left\{\int_t^\infty \left(u^{-\theta_1}a_1(u)K(u,f)\right)^{q_1}\frac{du}{u}\right\}^{\frac{1}{q_1}} \leq$$
$$\leq \left(\sup_{t<u<\infty} u^{-\theta_1}a_1(u)K(u,f)\right)^{\frac{q_1-r}{q_1}}\left\{\int_t^\infty \left(u^{-\theta_1}a_1(u)K(u,f)\right)^r\frac{du}{u}\right\}^{\frac{1}{q_1}} \prec$$
$$\prec \left(\left\|u^{-\theta_1-1/r}a_1(u)K(u,f)\right\|_{r,(t,\infty)}\right)^{\frac{q_1-r}{q_1}}\left(\left\|u^{-\theta_1-1/r}a_1(u)K(u,f)\right\|_{r,(t,\infty)}\right)^{\frac{r}{q_1}} =$$
$$= \left\|u^{-\theta_1-1/r}a_1(u)K(u,f)\right\|_{r,(t,\infty)}.$$

Using now (10) and applying Fubini's theorem we conclude $I_3 \prec I_1$. □

**Theorem 16.** Let $0 < \theta, \theta_1 < 1$, $0 < r, r_0, r_1, q_1 \leq \infty$, $a, a_1, b_0, b_1 \in SV$, $\left\|s^{-1/r_0}b_0(s)\right\|_{r_0,(1,\infty)} < \infty$, and $\left\|s^{-1/r_1}\frac{b_1(s)}{a_1(s)}\right\|_{r_1,(1,\infty)} < \infty$. Put $c_0(t) = \left\|s^{-1/r_0}b_0(s)\right\|_{r_0,(t,\infty)}$, $c_1(t) = a_1(t)\left\|s^{-1/r_1}\frac{b_1(s)}{a_1(s)}\right\|_{r_1,(t,\infty)}$, and $\rho(t) = t^{\theta_1}\frac{c_0(t)}{c_1(t)}$. Then
$$\left(\bar{A}_{0,r_0,b_0}, \bar{A}^{\mathcal{L}}_{\theta_1,r_1,b_1,q_1,a_1}\right)_{\theta,r,a} \cong \bar{A}_{\eta,r,a^{\#}},$$
where $\eta = \theta\theta_1$ and $a^{\#}(t) = c_0(t)^{1-\theta}c_1(t)^{\theta}a(\rho(t))$.

*Proof.* Denote $X_0 = A_0$, $X_1 = \bar{A}^{\mathcal{L}}_{\theta_1,r_1,b_1,q_1,a_1}$, $Z = \left(\bar{A}_{0,r_0,b_0}, \bar{A}^{\mathcal{L}}_{\theta_1,r_1,b_1,q_1,a_1}\right)_{\theta,r,a}$. Let $\sigma$ be a strongly increasing, differentiable function such that $\sigma(t) \approx t^{\theta_1}\frac{1}{c_1(t)}$. By Theorem 15 we have
$$\bar{A}_{0,r_0,b_0} \cong \left(A_0, \bar{A}^{\mathcal{L}}_{\theta_1,r_1,b_1,q_1,a_1}\right)_{0,r_0,d},$$
where $d = b_0 \circ \sigma^{(-1)}$. Hence, $Z \cong \left(\bar{X}_{0,r_0,d}, X_1\right)_{\theta,r,a}$. Using now Theorem 8 (i), we get $Z \cong \bar{X}_{\theta,r,c}$, where $c(t) = \left(\left\|s^{-1/r_0}d(s)\right\|_{r_0,(t,\infty)}\right)^{1-\theta}a\left(t\left\|s^{-1/r_0}d(s)\right\|_{r_0,(t,\infty)}\right)$. By the change of variables $s = \sigma(u)$, we arrive at
$$\left\|s^{-1/r_0}d(s)\right\|_{r_0,(\sigma(t),\infty)} \approx \left\|u^{-1/r_0}d(\sigma(u))\right\|_{r_0,(t,\infty)} = \left\|u^{-1/r_0}b_0(u)\right\|_{r_0,(t,\infty)} = c_0(t).$$
Because $c_1(t)^{\theta}c(\sigma(t)) \approx c_1(t)^{\theta}c_0(t)^{1-\theta}a(\sigma(t)c_0(t)) \approx a^{\#}(t)$, using once more Theorem 15, we obtain
$$Z \cong \left(A_0, \bar{A}^{\mathcal{L}}_{\theta_1,r_1,b_1,q_1,a_1}\right)_{\theta,r,c} \cong \bar{A}_{\eta,r,a^{\#}}.$$
This completes the proof. $\square$

**Theorem 17.** Let $0 < \theta_0 < \theta_1 < 1$, $0 < q_1, r_0, r_1 \leq \infty$, $a_1, b_0, b_1 \in SV$, $\left\|s^{-1/r_0}\frac{b_1(s)}{a_1(s)}\right\|_{r_0,(1,\infty)} < \infty$. Put $c_1(t) = a_1(t)\left\|s^{-1/r_1}\frac{b_1(s)}{a_1(s)}\right\|_{r_1,(t,\infty)}$ and $\rho(t) = t^{\theta_1-\theta_0}\frac{b_0(t)}{c_1(t)}$. Then for all $f \in \bar{A}_{\theta_0,r_0,b_0} + \bar{A}^{\mathcal{L}}_{\theta_1,r_1,b_1,q_1,a_1}$ and all $t > 0$,
$$K\left(\rho(t), f; \bar{A}_{\theta_0,r_0,b_0}, \bar{A}^{\mathcal{L}}_{\theta_1,r_1,b_1,q_1,a_1}\right) \approx \left\|u^{-\theta_0-1/r_0}b_0(u)K(u,f)\right\|_{r_0,(0,t)} +$$
$$+\rho(t)\left\|s^{-1/r_1}\frac{b_1(s)}{a_1(s)}\left\|u^{-\theta_1-1/q_1}a_1(u)K(u,f)\right\|_{q_1,(t,s)}\right\|_{r_1,(t,\infty)}. \tag{11}$$

*Proof.* Recall, we adopted the notation from the paper [9]. Let $\Phi_0$ be the function space corresponding to $\bar{A}_{\theta_0,r_0,b_0}$:
$$\|F(*)\|_{\Phi_0} = \left\|u^{-\theta_0-1/r_0}b_0(u)F(u)\right\|_{r_0,(0,\infty)}.$$
Let us calculate the functions $I(t,f)$, $g_0(t)$, and $h_0(t)$.
$$I(t,f) := \left\|\chi_{(0,t)}(*)K(*,f)\right\|_{\Phi_0} = \left\|u^{-\theta_0-1/r_0}b_0(u)K(u,f)\right\|_{r_0,(0,t)}.$$
$$g_0(t) := t\left\|\chi_{(t,\infty)}(*)\right\|_{\Phi_0} = t\left\|u^{-\theta_0-1/r_0}b_0(u)\right\|_{r_0,(t,\infty)} \approx t^{1-\theta_0}b_0(t).$$
$$h_0(t) := \left\|* \chi_{(0,t)}(*)\right\|_{\Phi_0} = \left\|u^{1-\theta_0-1/r_0}b_0(u)\right\|_{r_0,(0,t)} \approx t^{1-\theta_0}b_0(t).$$
Thus, $g_0(t) + h_0(t) \approx h_0(t) \approx g_0(t) \approx t^{1-\theta_0}b_0(t)$.

Let $\Phi_1$ be the function space corresponding to $\bar{A}^{\mathcal{L}}_{\theta_1,r_1,b_1,q_1,a_1}$:
$$\|F(*)\|_{\Phi_1} = \left\|s^{-1/r_1}\frac{b_1(s)}{a_1(s)}\left\|u^{-\theta_1-1/q_1}a_1(u)F(u)\right\|_{q_1,(0,s)}\right\|_{r_1,(0,\infty)}.$$
From the proof of Theorem 14 we have:
$$J(t,f) = \left\|s^{-1/r_1}\frac{b_1(s)}{a_1(s)}\left\|u^{-\theta_1-1/q_1}a_1(u)K(u,f)\right\|_{q_1,(t,s)}\right\|_{r_1,(t,\infty)},$$
$$g_1(t) \prec t^{1-\theta_1}c_1(t),$$
$$h_1(t) \approx t^{1-\theta_1}c_1(t).$$

Now we estimate the functions $g(t)$ and $h(t)$.

$$\frac{1}{g(t)} := \left\|\chi_{(t,\infty)}(u)\frac{u}{g_0(u)+h_0(u)}\right\|_{\Phi_1} \approx$$

$$\approx \left\|s^{-1/r_1}\frac{b_1(s)}{a_1(s)}\left\|u^{\theta_0-\theta_1-1/q_1}\frac{a_1(u)}{b_0(u)}\right\|_{q_1,(t,s)}\right\|_{r_1,(t,\infty)} \leq$$

$$\leq \left\|s^{-1/r_1}\frac{b_1(s)}{a_1(s)}\left\|u^{\theta_0-\theta_1-1/q_1}\frac{a_1(u)}{b_0(u)}\right\|_{q_1,(t,\infty)}\right\|_{r_1,(t,\infty)} \approx$$

$$\approx t^{\theta_0-\theta_1}\frac{a_1(t)}{b_0(t)}\left\|s^{-1/r_1}\frac{b_1(s)}{a_1(s)}\right\|_{r_1,(t,\infty)} = t^{\theta_0-\theta_1}\frac{c_1(t)}{b_0(t)}.$$

$$h(t) := \left\|\chi_{(0,t)}(u)\frac{u}{g_1(u)+h_1(u)}\right\|_{\Phi_0} \approx$$

$$\approx \left\|u^{\theta_1-\theta_0-1/r_0}\frac{b_0(u)}{c_1(u)}\right\|_{r_0,(0,t)} \approx t^{\theta_1-\theta_0}\frac{b_0(t)}{c_1(t)} = \rho(t).$$

We see that $h_0(t) \approx g_0(t)$, $h(t) \approx \rho(t) \prec g(t)$, and $\frac{h_0(t)}{h_1(t)} \approx \rho(t)$. Hence, by [9, Theorem 4, Case 4] we get (11). □

**Theorem 18.** Let $0 < \theta_0 < \theta_1 < 1$, $0 < q_1, r, r_1, r_0 \leq \infty$, $a, a_1, b_0, b_1 \in SV$, and $\left\|s^{-1/r_1}\frac{b_1(s)}{a_1(s)}\right\|_{r_1,(1,\infty)} < \infty$. Put $c_1(t) = a_1(t)\left\|s^{-1/r_1}\frac{b_1(s)}{a_1(s)}\right\|_{r_1,(t,\infty)}$ and $\rho(t) = t^{\theta_1-\theta_0}\frac{b_0(t)}{c_1(t)}$.

(i) If $0 < \theta < 1$ then
$$\left(\bar{A}_{\theta_0,r_0,b_0}, \bar{A}^{\mathcal{L}}_{\theta_1,r_1,b_1,q_1,a_1}\right)_{\theta,r,a} \cong \bar{A}_{\eta,r,a^{\#}},$$
where $\eta = (1-\theta)\theta_0 + \theta\theta_1$ and $a^{\#}(t) = b_0(t)^{1-\theta}c_1(t)^{\theta}a(\rho(t))$.

(ii) If additionally $\left\|s^{-1/r}a(s)\right\|_{r,(1,\infty)} < \infty$, then
$$\left(\bar{A}_{\theta_0,r_0,b_0}, \bar{A}^{\mathcal{L}}_{\theta_1,r_1,b_1,q_1,a_1}\right)_{0,r,a} \cong \bar{A}^{\mathcal{L}}_{\theta_0,r,a^{\#},r_0,b_0},$$
where $a^{\#}(t) = b_0(t)a(\rho(t))$.

*Proof.* Denote $X_0 = \bar{A}_{\theta_0,r_0,b_0}$, $X_1 = \bar{A}^{\mathcal{L}}_{\theta_1,r_1,b_1,q_1,a_1}$, $\bar{K}(t,f) = K(t,f;\bar{X})$, and $Z = \left(\bar{A}_{\theta_0,r_0,b_0}, \bar{A}^{\mathcal{L}}_{\theta_1,r_1,b_1,q_1,a_1}\right)_{\theta,r,a}$. We write for shortness $e(t) = b_0(t)^{-\theta}c_1(t)^{\theta}a(\rho(t))$. By Remark 3 (v) and Theorem 17 we can write

$$\|f\|_Z = \left\|u^{-\theta-1/r}a(u)\bar{K}(u,f)\right\|_{r,(0,\infty)} \approx$$

$$\approx \left\|\rho(t)^{-\theta}t^{-1/r}a(\rho(t))\bar{K}(\rho(t),f)\right\|_{r,(0,\infty)} =$$

$$= \left\|t^{-\theta(\theta_1-\theta_0)-1/r}e(t)\bar{K}(\rho(t),f)\right\|_{r,(0,\infty)} \approx I_1 + I_2,$$

where

$$I_1 := \left\|t^{-\theta(\theta_1-\theta_0)-1/r}e(t)\left\|u^{-\theta_0-1/r_0}b_0(u)K(u,f)\right\|_{r_0,(0,t)}\right\|_{r,(0,\infty)},$$

and

$$I_2 := \left\|t^{-\theta(\theta_1-\theta_0)-1/r}e(t)\rho(t)\left\|s^{-1/r_1}\frac{b_1(s)}{a_1(s)}\left\|u^{-\theta_1-1/q_1}a_1(u)K(u,f)\right\|_{q_1,(t,s)}\right\|_{r_1,(t,\infty)}\right\|_{r,(0,\infty)}.$$

Let $I_3 := \left\|t^{-\theta(\theta_1-\theta_0)-\theta_0-1/r}e(t)b_0(t)K(t,f)\right\|_{r,(0,\infty)}$. Because $t^{-1}K(t,f)\downarrow$ and $1-\theta_0 > 0$, we conclude
$$t^{-\theta_0}b_0(t)K(t,f) \approx t^{-1}K(t,f)\left\|u^{1-\theta_0-1/r_0}b_0(u)\right\|_{r_0,(0,t)} \leq \left\|u^{-\theta_0-1/r_0}b_0(u)K(u,f)\right\|_{r_0,(0,t)}.$$
Thus, $I_3 \prec I_1$.

(i) Let $0 < \theta < 1$. We have

$$I_3 = \left\| t^{-\eta-1/r} b_0(t)^{1-\theta} c_1(t)^\theta a(\rho(t)) K(t,f) \right\|_{r,(0,\infty)} = \|f\|_{\eta,r,a^\#}.$$

Thus, it is enough to show that $I_2 \prec I_3$ and $I_1 \prec I_3$. Consider $I_2$. Note that

$$\left\| s^{-1/r_1} \frac{b_1(s)}{a_1(s)} \left\| u^{-\theta_1-1/q_1} a_1(u) K(u,f) \right\|_{q_1,(t,s)} \right\|_{r_1,(t,\infty)} \leq$$
$$\left\| s^{-1/r_1} \frac{b_1(s)}{a_1(s)} \right\|_{r_1,(t,\infty)} \left\| u^{-\theta_1-1/q_1} a_1(u) K(u,f) \right\|_{q_1,(t,\infty)}.$$

Denote

$$d(t) = e(t) \frac{b_0(t)}{c_1(t)} \left\| s^{-1/r_1} \frac{b_1(s)}{a_1(s)} \right\|_{r_1,(t,\infty)} = b_0(t)^{1-\theta} c_1(t)^{\theta-1} \left\| s^{-1/r_1} \frac{b_1(s)}{a_1(s)} \right\|_{r_1,(t,\infty)} a(\rho(t)).$$

We get

$$I_2 \leq \left\| t^{(1-\theta)(\theta_1-\theta_0)-1/r} d(t) \left\| u^{-\theta_1-1/q_1} a_1(u) K(u,f) \right\|_{q_1,(t,\infty)} \right\|_{r,(0,\infty)} =$$
$$= \left\| t^{\theta_1-\eta-1/r} d(t) \left\| u^{-\theta_1-1/q_1} a_1(u) K(u,f) \right\|_{q_1,(t,\infty)} \right\|_{r,(0,\infty)} := I_4.$$

The estimates $I_4 \prec I_3$ and $I_1 \prec I_3$ can be proved by repeating the corresponding arguments of the proof of Theorem 15, using Lemma 4 (ii) or Lemma 4 (i).

(ii) Let now $\theta = 1$. In this case we have

$$I_1 = \left\| t^{-1/r} a(\rho(t)) \left\| u^{-\theta_0-1/r_0} b_0(u) K(u,f) \right\|_{r_0,(0,t)} \right\|_{r,(0,\infty)} = \|f\|_{\mathcal{L};\theta_0,r,a^\#,r_0,b_0},$$

$$I_2 = \left\| t^{\theta_1-\theta_0-1/r} a(\rho(t)) \frac{b_0(t)}{c_1(t)} \left\| s^{-1/r_1} \frac{b_1(s)}{a_1(s)} \left\| u^{-\theta_1-1/q_1} a_1(u) K(u,f) \right\|_{q_1,(t,s)} \right\|_{r_1,(t,\infty)} \right\|_{r,(0,\infty)},$$

and

$$I_3 = \left\| t^{-\theta_0-1/r} a(\rho(t)) b_0(t) K(t,f) \right\|_{r,(0,\infty)}.$$

Because $I_3 \prec I_1$ it is enough to show that $I_2 \prec I_3$. Denoting $d(t) = a(\rho(t)) \frac{b_0(t)}{a_1(t)}$, we get

$$I_2 \leq \left\| t^{(\theta_1-\theta_0)-1/r} a(\rho(t)) \frac{b_0(t)}{c_1(t)} \right.$$
$$\left. * \left\| s^{-1/r_1} \frac{b_1(s)}{a_1(s)} \right\|_{r_1,(t,\infty)} \left\| u^{-\theta_1-1/q_1} a_1(u) K(u,f) \right\|_{q_1,(t,\infty)} \right\|_{r,(0,\infty)} =$$
$$= \left\| t^{(\theta_1-\theta_0)-1/r} d(t) \left\| u^{-\theta_1-1/q_1} a_1(u) K(u,f) \right\|_{q_1,(t,\infty)} \right\|_{r,(0,\infty)} := I_4.$$

The estimate $I_4 \prec I_3$ can be proved as above by applying Hardy-type inequality and Fubini's theorem. □

## 5. Interpolation between the $\mathcal{R}$-spaces and the standard interpolation spaces

The proofs of all theorems in this section are based on the corresponding reiteration theorems from the previous section and on the formulae (6) and (7).

**Theorem 19.** Let $0 < \theta_0 < \theta_1 < 1$, $0 < q_1, r_0, r_1, r \leq \infty$, $a, a_1, b_0, b_1 \in SV$, and $\left\| s^{-1/r_1} \frac{b_1(s)}{a_1(s)} \right\|_{r_1,(0,1)} < \infty$. Put $c_1(t) = a_1(t) \left\| s^{-1/r_1} \frac{b_1(s)}{a_1(s)} \right\|_{r_1,(0,t)}$ and $\rho(t) = t^{\theta_1-\theta_0} \frac{b_0(t)}{c_1(t)}$.

(i) If $0 < \theta < 1$ then

$$\left( \bar{A}_{\theta_0,r_0,b_0}, \bar{A}^{\mathcal{R}}_{\theta_1,r_1,b_1,q_1,a_1} \right)_{\theta,r,a} \cong \bar{A}_{\eta,r,a^\#},$$

where $\eta = (1-\theta)\theta_0 + \theta\theta_1$ and $a^\#(t) \approx b_0(t)^{1-\theta} c_1(t)^\theta a(\rho(t))$.

(ii) If additionally $\left\| s^{-1/r} a(s) \right\|_{r,(1,\infty)} < \infty$, then

$$\left(\bar{A}_{\theta_0,r_0,b_0}, \bar{A}^{\mathcal{R}}_{\theta_1,r_1,b_1,q_1,a_1}\right)_{0,r,a} \cong \vec{A}^{\mathcal{L}}_{\theta_0,r,a^\#,r_0,b_0},$$

where $a^\#(t) \approx b_0(t)a(\rho(t))$.

**Theorem 20.** Let $0<\theta, \theta_1<1$, $0<q_1, r_0, r_1, r \leq \infty$, $a, a_1, b_0, b_1 \in SV$, $\left\|s^{-1/r_1}\frac{b_1(s)}{a_1(s)}\right\|_{r_1,(0,1)} < \infty$, and $\left\|s^{-1/r_0}b_0(s)\right\|_{r_0,(1,\infty)} < \infty$. Put $c_0(t) = \left\|s^{-1/r_0}b_0(s)\right\|_{r_0,(t,\infty)}$, $c_1(t) = a_1(t)\left\|s^{-1/r_1}\frac{b_1(s)}{a_1(s)}\right\|_{r_1,(0,t)}$, and $\rho(t) = t^{\theta_1}\frac{c_0(t)}{c_1(t)}$. Then
$$\left(\bar{A}_{0,r_0,b_0}, \bar{A}^{\mathcal{R}}_{\theta_1,r_1,b_1,q_1,a_1}\right)_{\theta,r,a} \cong \bar{A}_{\eta,r,a^\#},$$
where $\eta = \theta\theta_1$ and $a^\#(t) = c_0(t)^{1-\theta}c_1(t)^\theta a(\rho(t))$.

**Theorem 21.** Let $0 < \theta_1 < 1$, $0 \leq \theta < 1$, $0 < q_1, r_1, r \leq \infty$, $a, a_1, b_1 \in SV$, and $\left\|s^{-1/r_1}\frac{b_1(s)}{a_1(s)}\right\|_{r_1,(0,1)} < \infty$. If $\theta = 0$, we additionally assume that $\left\|s^{-1/r}a(s)\right\|_{r,(1,\infty)} < \infty$. Put $c_1(t) = a_1(t)\left\|s^{-1/r_1}\frac{b_1(s)}{a_1(s)}\right\|_{r_1,(0,t)}$ and $(t) = t^{\theta_1}\frac{1}{c_1(t)}$. Then
$$\left(\bar{A}_0, \bar{A}^{\mathcal{R}}_{\theta_1,r_1,b_1,q_1,a_1}\right)_{\theta,r,a} \cong \bar{A}_{\eta,r,a^\#},$$
where $\eta = \theta\theta_1$ and $a^\#(t) = c_1(t)^\theta a(\rho(t))$.

**Theorem 22.** Let $0 < \theta_0 < 1$, $0 < \theta \leq 1$, $0 < q_0, r, r_0 \leq \infty$, $a, a_0, b_0 \in SV$, and $\left\|s^{-1/r_0}\frac{b_0(s)}{a_0(s)}\right\|_{r_0,(0,1)} < \infty$. If $\theta = 1$, we additionally assume that $\left\|s^{-1/r}a(s)\right\|_{r,(0,1)} < \infty$. Put $c_0(t) = a_0(t)\left\|s^{-1/r_0}\frac{b_0(s)}{a_0(s)}\right\|_{r_1,(0,t)}$ and $\rho(t) = t^{1-\theta_0}c_0(t)$. Then
$$\left(\bar{A}^{\mathcal{R}}_{\theta_0,r_0,b_0,q_0,a_0}, A_1\right)_{\theta,r,a} \cong \bar{A}_{\eta,r,a^\#},$$
where $\eta = (1-\theta)\theta_0 + \theta$ and $a^\#(t) = c_0(t)^{1-\theta}a(\rho(t))$.

**Theorem 23.** Let $0 < \theta, \theta_0 < 1$, $0 < q_0, r, r_0, r_1 \leq \infty$, $a, b_0, b_1, a_0 \in SV$, $\left\|s^{-1/r_0}\frac{b_0(s)}{a_0(s)}\right\|_{r_0,(0,1)} < \infty$, and $\left\|s^{-1/r_1}b_1(s)\right\|_{r_1,(0,1)} < \infty$. Put $c_0(t) = a_0(t)\left\|s^{-1/r_0}\frac{b_0(s)}{a_0(s)}\right\|_{r_0,(0,t)}$, $c_1(t) = \left\|s^{-1/r_1}b_1(s)\right\|_{r_1,(0,t)}$, and $\rho(t) = t^{1-\theta_0}\frac{c_0(t)}{c_1(t)}$. Then
$$\left(\bar{A}^{\mathcal{R}}_{\theta_0,r_0,b_0,q_0,a_0}, \bar{A}_{1,r_1,b_1}\right)_{\theta,r,a} \cong \bar{A}_{\eta,r,a^\#},$$
where $\eta = (1-\theta)\theta_0 + \theta$ and $a^\#(t) = c_0(t)^{1-\theta}c_1(t)^\theta a(\rho(t))$.

**Theorem 24.** Let $0 < \theta_0 < \theta_1 < 1$, $0 < q_0, r, r_1, r_0 \leq \infty$, $a, a_0, b_0, b_1 \in SV$, and $\left\|s^{-1/r_0}\frac{b_0(s)}{a_0(s)}\right\|_{r_0,(0,1)} < \infty$. Put $c_0(t) = a_0(t)\left\|s^{-1/r_0}\frac{b_0(s)}{a_0(s)}\right\|_{r_0,(0,t)}$ and $\rho(t) = t^{\theta_1-\theta_0}\frac{c_0(t)}{b_1(t)}$.

(i) If $0 < \theta < 1$ then
$$\left(\bar{A}^{\mathcal{R}}_{\theta_0,r_0,b_0,q_0,a_0}, \bar{A}_{\theta_1,r_1,b_1}\right)_{\theta,r,a} \cong \bar{A}_{\eta,r,a^\#},$$
where $\eta = (1-\theta)\theta_0 + \theta\theta_1$ and $a^\#(t) = c_0(t)^{1-\theta}b_1(t)^\theta a(\rho(t))$.

(ii) If additionally $\left\|s^{-1/r}a(s)\right\|_{r,(0,1)} < \infty$, then
$$\left(\bar{A}^{\mathcal{R}}_{\theta_0,r_0,b_0,q_0,a_0}, \bar{A}_{\theta_1,r_1,b_1}\right)_{1,r,a} \cong \bar{A}^{\mathcal{R}}_{\theta_1,r,a^\#,r_1,b_1},$$
where $a^\#(t) = b_1(t)a(\rho(t))$.

## 6. Reiteration formulae for couples where both operands are $\mathcal{L}$ or $\mathcal{R}$ spaces

**Theorem 25.** Let $0 < \theta_0 < \theta_1 < 1$, $0 < \theta < 1$, $0 < r, r_0, r_1, q_0, q_1 \leq \infty$, $a, a_0, a_1, b_0, b_1 \in SV$, $\left\|s^{-1/r_0}\frac{b_0(s)}{a_0(s)}\right\|_{r_0,(1,\infty)} < \infty$, and $\left\|s^{-1/r_1}\frac{b_1(s)}{a_1(s)}\right\|_{r_1,(1,\infty)} < \infty$. Put $c_0(t) = a_0(t)\left\|s^{-1/r_0}\frac{b_0(s)}{a_0(s)}\right\|_{r_0,(t,\infty)}$, $c_1(t) = a_1(t)\left\|s^{-1/r_1}\frac{b_1(s)}{a_1(s)}\right\|_{r_1,(t,\infty)}$, and $\rho(t) = t^{\theta_1-\theta_0}\frac{c_0(t)}{c_1(t)}$. Then
$$\left(\bar{A}^{\mathcal{L}}_{\theta_0,r_0,b_0,q_0,a_0}, \bar{A}^{\mathcal{L}}_{\theta_1,r_1,b_1,q_1,a_1}\right)_{\theta,r,\mathrm{a}} \cong \bar{A}_{\eta,r,a^{\#}},$$
where $\eta = (1-\theta)\theta_0 + \theta\theta_1$ and $a^{\#}(t) = c_0(t)^{1-\theta}c_1(t)^{\theta}a(\rho(t))$.

*Proof.* Let $X_0 = \bar{A}_{\theta_0,q_0,a_0}$, $X_1 = \bar{A}^{\mathcal{L}}_{\theta_1,r_1,b_1,q_1,a_1}$, $Z = \left(\bar{A}^{\mathcal{L}}_{\theta_0,r_0,b_0,q_0,a_0}, \bar{A}^{\mathcal{L}}_{\theta_1,r_1,b_1,q_1,a_1}\right)_{\theta,r,\mathrm{a}}$. Let $\sigma$ be a strongly increasing, differentiable function such that $\sigma(t) \approx t^{\theta_1-\theta_0}\frac{a_0(t)}{c_1(t)}$ and $c(t) = \frac{b_0(\sigma^{(-1)}(t))}{a_0(\sigma^{(-1)}(t))}$. By Theorem 18 (ii) we have
$$\bar{A}^{\mathcal{L}}_{\theta_0,r_0,b_0,q_0,a_0} \cong \left(\bar{A}_{\theta_0,q_0,a_0}, \bar{A}^{\mathcal{L}}_{\theta_1,r_1,b_1,q_1,a_1}\right)_{0,r_0,c} = \bar{X}_{0,r_0,c}.$$
Hence, $Z \cong \left(\bar{X}_{0,r_0,c}, X_1\right)_{\theta,r,\mathrm{a}}$. Because $\left\|s^{-1/r_0}\frac{b_0(s)}{a_0(s)}\right\|_{r_0,(1,\infty)} < \infty$, it holds $\left\|s^{-1/r_0}c(t)\right\|_{r_0,(1,\infty)} < \infty$. This means that we can apply Theorem 8 (i). So, we get $Z \cong \bar{X}_{\theta,r,d}$, where $d(t) = \left(\left\|s^{-1/r_0}c(t)\right\|_{r_0,(t,\infty)}\right)^{1-\theta}a\left(t\left\|s^{-1/r_0}c(t)\right\|_{r_0,(t,\infty)}\right)$. By Theorem 18 (i) we obtain
$$Z \cong \left(\bar{A}_{\theta_0,q_0,a_0}, \bar{A}^{\mathcal{L}}_{\theta_1,r_1,b_1,q_1,a_1}\right)_{\theta,r,d} \cong \bar{A}_{\eta,r,a^{\#}},$$
where $a^{\#}(t) = a_0(t)^{1-\theta}c_1(t)^{\theta}d(\sigma(t))$. By the change of variables $u = \sigma^{(-1)}(s)$ we arrive at
$$\left\|s^{-1/r_0}c(s)\right\|_{r_0,(\sigma(t),\infty)} = \left\|s^{-1/r_0}\frac{b_0(\sigma^{(-1)}(s))}{a_0(\sigma^{(-1)}(s))}\right\|_{r_0,(\sigma(t),\infty)} \approx \left\|u^{-1/r_0}\frac{b_0(u)}{a_0(u)}\right\|_{r_0,(t,\infty)}.$$
Thus, because $\sigma(t)\left\|s^{-1/r_0}\frac{b_0(s)}{a_0(s)}\right\|_{r_0,(t,\infty)} \approx \rho(t)$, we conclude
$$d(\sigma(t)) = \left(\left\|s^{-1/r_0}c(s)\right\|_{r_0,(\sigma(t),\infty)}\right)^{1-\theta} a\left(\sigma(t)\left\|s^{-1/r_0}c(s)\right\|_{r_0,(\sigma(t),\infty)}\right) \approx$$
$$\approx \left(\left\|s^{-1/r_0}\frac{b_0(s)}{a_0(s)}\right\|_{r_0,(t,\infty)}\right)^{1-\theta} a(\rho(t)).$$
Hence, $a^{\#}(t) \approx c_0(t)^{1-\theta}c_1(t)^{\theta}a(\rho(t))$. □

Repeating the previous proof with Theorem 19 instead of Theorem 18 leads us to the following statement.

**Theorem 26.** Let $0 < \theta_0 < \theta_1 < 1$, $0 < \theta < 1$, $0 < r, r_0, r_1, q_0, q_1 \leq \infty$, $a, a_0, a_1, b_0, b_1 \in SV$, $\left\|s^{-1/r_0}\frac{b_0(s)}{a_0(s)}\right\|_{r_0,(1,\infty)} < \infty$, and $\left\|s^{-1/r_1}\frac{b_1(s)}{a_1(s)}\right\|_{r_1,(0,1)} < \infty$. Put $c_0(t) = a_0(t)\left\|s^{-1/r_0}\frac{b_0(s)}{a_0(s)}\right\|_{r_0,(t,\infty)}$, $c_1(t) = a_1(t)\left\|s^{-1/r_1}\frac{b_1(s)}{a_1(s)}\right\|_{r_1,(0,t)}$, and $\rho(t) = t^{\theta_1-\theta_0}\frac{c_0(t)}{c_1(t)}$. Then
$$\left(\bar{A}^{\mathcal{L}}_{\theta_0,r_0,b_0,q_0,a_0}, \bar{A}^{\mathcal{R}}_{\theta_1,r_1,b_1,q_1,a_1}\right)_{\theta,r,\mathrm{a}} \cong \bar{A}_{\eta,r,a^{\#}},$$

Where $\eta = (1-\theta)\theta_0 + \theta\theta_1$ and $a^\#(t) = c_0(t)^{1-\theta} c_1(t)^\theta a(\rho(t))$.

Using the symmetry arguments given by the formulae (6) and (7), the following theorem can be proved.

**Theorem 27.** Let $0 < \theta_0 < \theta_1 < 1$, $0 < \theta < 1$, $0 < r, r_0, r_1, q_0, q_1 \leq \infty$, $a, a_0, a_1, b_0, b_1 \in SV$, $\left\| s^{-1/r_0} \frac{b_0(s)}{a_0(s)} \right\|_{r_0,(0,1)} < \infty$, and $\left\| s^{-1/r_1} \frac{b_1(s)}{a_1(s)} \right\|_{r_1,(0,1)} < \infty$. Put $c_0(t) = a_0(t) \left\| s^{-1/r_0} \frac{b_0(s)}{a_0(s)} \right\|_{r_0,(0,t)}$, $c_1(t) = a_1(t) \left\| s^{-1/r_1} \frac{b_1(s)}{a_1(s)} \right\|_{r_1,(0,t)}$, and $\rho(t) = t^{\theta_1 - \theta_0} \frac{c_0(t)}{c_1(t)}$. Then
$$\left( \bar{A}^{\mathcal{R}}_{\theta_0,r_0,b_0,q_0,a_0}, \bar{A}^{\mathcal{R}}_{\theta_1,r_1,b_1,q_1,a_1} \right)_{\theta,r,a} \cong \bar{A}_{\eta,r,a^\#},$$
where $\eta = (1-\theta)\theta_0 + \theta\theta_1$ and $a^\#(t) = c_0(t)^{1-\theta} c_1(t)^\theta a(\rho(t))$.

The following theorem can be proved can be proved directly using relevant Holmstedt-type formula for the $K$-functional.

**Theorem 28.** Let $0 < \theta_0 < \theta_1 < 1$, $0 < \theta < 1$, $0 < r, r_0, r_1, q_0, q_1 \leq \infty$, $a, a_0, a_1, b_0, b_1 \in SV$, $\left\| s^{-1/r_0} \frac{b_0(s)}{a_0(s)} \right\|_{r_0,(0,1)} < \infty$, and $\left\| s^{-1/r_1} \frac{b_1(s)}{a_1(s)} \right\|_{r_1,(1,\infty)} < \infty$. Put $c_0(t) = a_0(t) \left\| s^{-1/r_0} \frac{b_0(s)}{a_0(s)} \right\|_{r_0,(0,t)}$, $c_1(t) = a_1(t) \left\| s^{-1/r_1} \frac{b_1(s)}{a_1(s)} \right\|_{r_1,(t,\infty)}$, and $\rho(t) = t^{\theta_1 - \theta_0} \frac{c_0(t)}{c_1(t)}$. Then
$$\left( \bar{A}^{\mathcal{R}}_{\theta_0,r_0,b_0,q_0,a_0}, \bar{A}^{\mathcal{L}}_{\theta_1,r_1,b_1,q_1,a_1} \right)_{\theta,r,a} \cong \bar{A}_{\eta,r,a^\#},$$
where $\eta = (1-\theta)\theta_0 + \theta\theta_1$ and $a^\#(t) = c_0(t)^{1-\theta} c_1(t)^\theta a(\rho(t))$.

## 7. Ordered couples

In this section, we assume that the couple $(A_0, A_1)$ is ordered in the sense that $A_1 \subset A_0$. This is the case for example if $A_0 = L_1(\Omega)$ and $A_1 = L_\infty(\Omega)$ where $\Omega \subset R^n$ with $|\Omega| < \infty$. We briefly outline how the previous results need to be amended for this particular case.

Since $A_1 \subset A_0$, it holds
$$K(u, f) \approx \|f_0\|_{A_0}, \quad u > 1. \tag{12}$$
Based on (12), it is natural to define interpolation spaces with integration over $(0,1)$. (See the arguments in [6, Section 7].)

**Definition 29.** Let $A_1 \subset A_0$, $0 \leq \theta \leq 1$, $0 < q \leq \infty$ and $b \in SV(0,1)$. We put
$$\bar{A}^{(0,1)}_{\theta,q;b} \equiv (A_0, A_1)^{(0,1)}_{\theta,q;b} :=$$
$$:= \left\{ f \in A_0 + A_1 : \|f\|_{\theta,q;b;(0,1)} = \left\| u^{-\theta - 1/q} b(u) K(u, f) \right\|_{q,(0,1)} < \infty \right\}.$$

**Definition 30.** Let $A_1 \subset A_0$, $0 < r, q \leq \infty$, $0 < \sigma < 1$, $a, b \in SV(0,1)$. We put
$$\bar{A}^{\mathcal{L};(0,1)}_{\sigma,r,b,q,a} := \left\{ f \in A_0 + A_1 : \|f\|_{\mathcal{L};\sigma,r,b,q,a;(0,1)} := \right.$$
$$\left. \left\| t^{-1/r} \frac{b(t)}{a(t)} \left\| u^{-\sigma - 1/q} a(u) K(u, f) \right\|_{q,(0,t)} \right\|_{r,(0,1)} < \infty \right\}.$$

Similarly,
$$\bar{A}^{\mathcal{R};(0,1)}_{\sigma,r,b,q,a} := \left\{ f \in A_0 + A_1 : \|f\|_{\mathcal{R};\sigma,r,b,q,a;(0,1)} := \right.$$
$$\left. \left\| t^{-1/r} \frac{b(t)}{a(t)} \left\| u^{-\sigma - 1/q} a(u) K(u, f) \right\|_{q,(t,1)} \right\|_{r,(0,1)} < \infty \right\}.$$

Using (12), it is not difficult to show that

(i)   $\bar{A}_{\theta,q;b} \cong \bar{A}_{\theta,q;b}^{(0,1)}$,

(ii)  if $\left\|s^{-1/q}b(s)\right\|_{q,(0,1)} < \infty$, then $\bar{A}_{0,q;b}^{(0,1)} \cong A_0$,

(iii) $\bar{A}_{\sigma,r,b,q,a}^{\mathcal{L};(0,1)} \subset \bar{A}_{\sigma,q;a}^{(0,1)}$,

(iv)  if $\left\|s^{-1/r}\frac{b(s)}{a(s)}\right\|_{r,(0,1)} < \infty$, then $\bar{A}_{\sigma,r,b,q,a}^{\mathcal{L};(0,1)} \cong \bar{A}_{\sigma,q;a}^{(0,1)}$,

(v)   if $\left\|s^{-1/r}\frac{b(s)}{a(s)}\right\|_{r,(0,1)} = \infty$, then $\bar{A}_{\sigma,r,b,q,a}^{\mathcal{R};(0,1)} = \{0\}$,

(vi)  if $\left\|s^{-1/r}\frac{b(s)}{a(s)}\right\|_{r,(0,1)} < \infty$, then $\bar{A}_{\sigma,q;b}^{(0,1)} \subset \bar{A}_{\sigma,r,b,q,a}^{\mathcal{R};(0,1)}$,

(vii) if $\left\|s^{-1/r}\frac{b(s)}{a(s)}\right\|_{r,(1,\infty)} < \infty$, then $\bar{A}_{\sigma,r,b,q,a}^{\mathcal{L}} \cong \bar{A}_{\sigma,r,b,q,a}^{\mathcal{L};(0,1)}$,

(viii) if $\left\|s^{-1/r}\frac{b(s)}{a(s)}\right\|_{r,(0,1)} < \infty$, then $\bar{A}_{\sigma,r,b,q,a}^{\mathcal{R}} \cong \bar{A}_{\sigma,r,b,q,a}^{\mathcal{R};(0,1)}$.

In the properties (i), (vii), and (viii) we assume that the functions from $SV(0,1)$ are extended to suitable functions from $SV(0,\infty)$.

These properties show how all main theorems from previous sections can be reformulated for the case when $A_1 \subset A_0$. For example, Theorem 28 can be reformulated as follows.

**Theorem 31.** Let $A_1 \subset A_0$, $0 < \theta_0 < \theta_1 < 1$, $0 < \theta < 1$, $0 < r, r_0, r_1, q_0, q_1 \leq \infty$, $a, a_0, a_1, b_0, b_1 \in SV(0,1)$, $\left\|s^{-1/r_0}\frac{b_0(s)}{a_0(s)}\right\|_{r_0,(0,1)} < \infty$, and $\left\|s^{-1/r_1}\frac{b_1(s)}{a_1(s)}\right\|_{r_1,(0,1)} = \infty$. Put
$c_0(t) = a_0(t)\left\|s^{-1/r_0}\frac{b_0(s)}{a_0(s)}\right\|_{r_0,(0,t)}$, $c_1(t) = a_1(t)\left(1 + \left\|s^{-1/r_1}\frac{b_1(s)}{a_1(s)}\right\|_{r_1,(t,1)}\right)$, and $\rho(t) = t^{\theta_1-\theta_0}\frac{c_0(t)}{c_1(t)}$ ($0 < t < 1$). Then
$$\left(\bar{A}_{\theta_0,r_0,b_0,q_0,a_0}^{\mathcal{R};(0,1)}, \bar{A}_{\theta_1,r_1,b_1,q_1,a_1}^{\mathcal{L};(0,1)}\right)_{\theta,r,a} \cong \bar{A}_{\eta,r,a^{\#}}^{(0,1)},$$
where $\eta = (1-\theta)\theta_0 + \theta\theta_1$ and $a^{\#}(t) = c_0(t)^{1-\theta}c_1(t)^{\theta}a(\rho(t))$.

## 8. Applications

Here we present some interpolation results for the grand and small Lebesgue spaces as applications of our general reiteration theorems. Let $(\Omega, \mu)$ denote a totally $\sigma$-finite measure space with a non-atomic measure $\mu$. We consider functions $f$ from the set $\mathfrak{M}(\Omega, \mu)$ of all $\mu$-measurable functions on $\Omega$. As conventional (see e.g. [1]), $f^*(t)$ ($t>0$) denotes the non-increasing rearrangement of $f$.

### 8.1. The Lorentz–Karamata spaces, $L^{\mathcal{L}}$ and $L^{\mathcal{R}}$ spaces, the grand and small Lorentz Spaces

**Definition 32.** ([6, Definition 5.1]). Let $0 < p,q \leq \infty$ and $b \in SV$. The Lorentz–Karamata space $L_{p,q;b} \equiv L_{p,q;b}(\Omega)$ is the set of all $f \in \mathfrak{M}(\Omega, \mu)$ such that
$$\|f\|_{p,q;a} := \left\|t^{\frac{1}{p}-\frac{1}{q}}b(t)f^*(t)\right\|_{q,(0,\infty)} < \infty.$$

The Lorentz–Karamata space $L_{p,q;1}(\Omega)$ coincides with the Lorentz space $L_{p,q}(\Omega)$; it becomes the Lebesgue space $L_p(\Omega)$ if $b = 1$ and $p = q$.

**Definition 33.** (Cf. [6, (5.21), (5.33)] [3, Definition 7.1].) Let $0<p<\infty$, $0<q,r\leq \infty$, and $a,b \in SV$. The spaces $L^{\mathcal{L}}_{p,r,b,q,a} \equiv L^{\mathcal{L}}_{p,r,b,q,a}(\Omega)$ and $L^{\mathcal{R}}_{p,r,b,q,a} \equiv L^{\mathcal{R}}_{p,r,b,q,a}(\Omega)$ are the sets of all $f \in \mathfrak{M}(\Omega,\mu)$ such that

$$\|f\|_{L^{\mathcal{L}}_{p,r,b,q,a}} := \left\| t^{-1/r} \frac{b(t)}{a(t)} \left\| u^{\frac{1}{p}-\frac{1}{q}} a(u) f^*(u) \right\|_{q,(0,t)} \right\|_{r,(0,\infty)} < \infty$$

or

$$\|f\|_{L^{\mathcal{R}}_{pr,b,q,a}} := \left\| t^{-1/r} \frac{b(t)}{a(t)} \left\| u^{\frac{1}{p}-\frac{1}{q}} a(u) f^{**}(u) \right\|_{q,(t,\infty)} \right\|_{r,(0,\infty)} < \infty,$$

correspondently.

We will require that $\left\| t^{-1/r} \frac{b(t)}{a(t)} \right\|_{r,(1,\infty)} < \infty$ for $L^{\mathcal{L}}_{p,r,b,q,a}$ and $\left\| t^{-1/r} \frac{b(t)}{a(t)} \right\|_{r,(01)} < \infty$ for $L^{\mathcal{R}}_{p,r,b,q,a}$. Otherwise the corresponding space consists only of the null-element. Similar definitions can be found in [4, 5, 8]. We refer to the spaces $L^{\mathcal{L}}_*$ and $L^{\mathcal{R}}_*$ as $L^{\mathcal{L}}$ and $L^{\mathcal{R}}$ spaces, respectively.

In order to be able to compare our results with the results from [8, 11, 13] we introduce the analogues of the grand and small Lorentz spaces. They find many different important applications and have been well studied by many authors, see [8, 11, 13] and the references therein. One can find different definitions of grand and small Lorentz spaces. Normally, they are defined for a bounded Lebesgue measurable domain $\Omega$ in $R^n$ with measure 1. In [13] it is required that the functions are real-valued. We do not require that $\mu(\Omega) = 1$ and that the functions are real-valued. Note, that using arguments of Section 8 all corollaries below can be easy reformulated for the case when $\mu(\Omega) = 1$.

**Definition 34.** Let $0<p<\infty$, $0<q,r\leq \infty$, and $b \in SV$. The small Lorentz space $L^{(p,q,r)}_b$ is the set of all $f \in \mathfrak{M}(\Omega,\mu)$ such that

$$\|f\|_{L^{(p,q,r)}_b(\Omega)} := \left\| t^{-\frac{1}{r}} b(t) \left\| u^{\frac{1}{p}-\frac{1}{q}} f^*(u) \right\|_{q,(0,t)} \right\|_{r,(0,\infty)} < \infty.$$

The grand Lorentz space $L^{p),q,r}_b$ is the set of all $f \in \mathfrak{M}(\Omega,\mu)$ such that

$$\|f\|_{L^{p),q,r}_b(\Omega)} := \left\| t^{-\frac{1}{r}} b(t) \left\| u^{\frac{1}{p}-\frac{1}{q}} f^{**}(u) \right\|_{q,(t,\infty)} \right\|_{r,(0,\infty)} < \infty.$$

It is clear, that $L^{(p,q,r)}_b = L^{\mathcal{L}}_{p,r,b,q,1}$ and $L^{p),q,r}_b = L^{\mathcal{R}}_{p,r,b,q,1}$.

For simplicity, we consider below only interpolation spaces between $L^1$ and $L^\infty$. Lemma 35 and Lemma 36 characterise the $L^{\mathcal{L}}$ and $L^{\mathcal{R}}$ spaces as appropriate $\mathcal{L}$- and $\mathcal{R}$- limiting interpolation spaces. Lemma 35 is a special case of [6, Lemma 5.4].

**Lemma 35.** Let $1 < p < \infty$, $0 < q,r \leq \infty$, $a,b \in SV$, and $\theta = 1 - \frac{1}{p}$. Then
$$L^{\mathcal{L}}_{p,r,b,q,a} \cong (L^1, L^\infty)^{\mathcal{L}}_{\theta,r,b,q,a}.$$
In particular, $L^{(p,q,r)}_b \cong (L^1, L^\infty)^{\mathcal{L}}_{\theta,r,b,q,1}$.

**Lemma 36.** ([6, Lemma 5.9].) Let $1 < p < \infty$, $0 < q,r \leq \infty$, $a,b \in SV$, and $\theta = 1 - \frac{1}{p}$. Then

$$L^{\mathcal{R}}_{p,r,b,q,a} \cong (L^1, L^\infty)^{\mathcal{R}}_{\theta,r,b,q,a},$$
In particular, $L^{p),q,r}_b \cong (L^1, L^\infty)^{\mathcal{R}}_{\theta,r,b,q,1}$.

Using Lemma 35 and Lemma 36, we are able to characterise interpolation spaces lying between $L^{\mathcal{L}}$ and $L^{\mathcal{R}}$ spaces and Lorentz–Karamata spaces. In corollaries below, we restrict ourselves to grand and small Lorentz spaces. Moreover, for $0 < \theta < 1$, we formulate the formulae only for the classical case $(*,*)_{\theta,r}$.

### 8.2. Interpolation between the small Lorentz space and the Lorentz spaces

**Corollary 37.** (Cf. [13, Theorem 5.3].) Let $0 < \theta < 1$, $1 < p_0 < p_1 < \infty$, $0 < r, r_0, r_1, q_0 \leq \infty$, $b_0 \in SV$, and $\|s^{-1/r_0} b_0(s)\|_{r_0,(1,\infty)} < \infty$. Put $c_0(t) = \|s^{-1/r_0} b_0(s)\|_{r_0,(t,\infty)}$. Then
$$\left( L^{(p_0,q_0,r_0)}_{b_0}, L_{p_1,r_1} \right)_{\theta,r} \cong L_{p,r,a^\#},$$
where $\frac{1}{p} = \frac{1-\theta}{p_0} + \frac{\theta}{p_1}$ and $a^\#(t) = c_0(t)^{1-\theta}$.

*Proof.* Let $\bar{A} = (L_1, L_\infty)$, $\theta_i = 1 - \frac{1}{p_i}$, and $\eta = (1-\theta)\theta_0 + \theta\theta_1$. By Lemma 35 we have
$$L^{(p_0,q_0,r_0)}_{b_0} \cong (L_1, L_\infty)^{\mathcal{L}}_{\theta_0,r_0,b_0,q_0,1} = \bar{A}^{\mathcal{L}}_{\theta_0,r_0,b_0,q_0,1}.$$
It is known [6, Corollary 5.3], that $\bar{A}_{\theta_1,r_1} \cong L_{p_1,r_1}$. Using now Theorem 11 (i) and [6, Corollary 5.3], we get
$$\left( L^{(p_0,q_0,r_0)}_{b_0}, L_{p_1,r_1} \right)_{\theta,r} \cong \left( \bar{A}^{\mathcal{L}}_{\theta_0,r_0,b_0,q_0,1}, \bar{A}_{\theta_1,r_1} \right)_{\theta,r} \cong \bar{A}_{\eta,r,a^\#} \cong L_{p,r,a^\#}.$$
This completes the proof. □

All corollaries below in this and the next subsections can be proved analogously using theorems from Sections 5 and 6.

**Corollary 38.** Let $1 < p_0 < p_1 < \infty$, $0 < r, r_0, r_1, q_0 \leq \infty$, $a, b_0 \in SV$, $\|s^{-1/r} a(s)\|_{r,(0,1)} < \infty$, and $\|s^{-1/r_0} b_0(s)\|_{r_0,(1,\infty)} < \infty$. Put $c_0(t) = \|s^{-1/r_0} b_0(s)\|_{r_0,(t,\infty)}$ and $\rho(t) = t^{\frac{1}{p_0} - \frac{1}{p_1}} c_0(t)$. Then
$$\left( L^{(p_0 q_0, r_0)}_{b_0}, L_{p_1,r_1} \right)_{1,r;a} \cong L^{p_1),r_1,r}_{a^\#},$$
where $a^\#(t) = a(\rho(t))$.

**Corollary 39.** (Cf. [8, Corollary 7.8].) Let $0 < \theta < 1$, $1 < p_0 < \infty$, $0 < r, r_0, q_0 \leq \infty$, and $\|s^{-1/r_0} b_0(s)\|_{r_0,(1,\infty)} < \infty$. Put $c_0(t) = \|s^{-1/r_0} b_0(s)\|_{r_0,(t,\infty)}$. Then
$$\left( L^{(p_0,q_0,r_0)}_{b_0}, L_\infty \right)_{\theta,r} \cong L_{p,r,a^\#},$$
where $\frac{1}{p} = \frac{1-\theta}{p_0}$ and $a^\#(t) = c_0(t)^{1-\theta}$.

**Corollary 40.** Let $0 < \theta < 1$, $1 < p_0 < p_1 < \infty$, $0 < r, r_0, r_1, q_1 \leq \infty$, $b_1 \in SV$, and $\|s^{-1/r_1} b_1(s)\|_{r_1,(1,\infty)} < \infty$. Put $c_1(t) = \|s^{-1/r_1} b_1(s)\|_{r_1,(t,\infty)}$. Then
$$\left( L_{p_0,r_0}, L^{(p_1,q_1,r_1)}_{b_1} \right)_{\theta,r} \cong L_{p,r,a^\#},$$
where $\frac{1}{p} = \frac{1-\theta}{p_0} + \frac{\theta}{p_1}$ and $a^\#(t) = c_1(t)^\theta$.

**Corollary 41.** Let $1 < p_0 < p_1 < \infty$, $0 < r, r_0, r_1, q_1 \leq \infty$, $a, b_1 \in SV$, $\|s^{-1/r_1} b_1(s)\|_{r_1,(1,\infty)} < \infty$, and $\|s^{-1/r} a(s)\|_{r,(1,\infty)} < \infty$. Put $c_1(t) = \|s^{-1/r_1} b_1(s)\|_{r_1,(t,\infty)}$ and $\rho(t) = t^{\frac{1}{p_0}-\frac{1}{p_1}} \frac{1}{c_1(t)}$. Then

$$\left(L_{p_0,r_0}, L_{b_1}^{(p_1,q_1,r_1)}\right)_{0,r;a} \cong L_{a^\#}^{(p_0,r_0,r)},$$

where $a^\#(t) = a(\rho(t))$.

**Corollary 42.** Let $0 < \theta < 1$, $1 < p_1 < \infty$, $0 < r, r_1, q_1 \leq \infty$, and $\|s^{-1/r_1} b_1(s)\|_{r_1,(1,\infty)} < \infty$. Put $c_1(t) = \|s^{-1/r_1} b_1(s)\|_{r_1,(t,\infty)}$. Then

$$\left(L_1, L_{b_1}^{(p_1,q_1,r_1)}\right)_{\theta,r} \cong L_{p,r,a^\#},$$

where $\frac{1}{p} = 1 - \theta + \frac{\theta}{p_1}$ and $a^\#(t) = c_1(t)^\theta$.

### 8.3. Interpolation between the grand Lorentz spaces and the Lorentz spaces

**Corollary 43.** Let $0 < \theta < 1$, $1 < p_0 < p_1 < \infty$, $0 < r, r_0, r_1, q_1 \leq \infty$, $b_1 \in SV$, and $\|s^{-1/r_1} b_1(s)\|_{r_1,(0,1)} < \infty$. Put $c_1(t) = \|s^{-1/r_1} b_1(s)\|_{r_1,(0,t)}$. Then

$$\left(L_{p_0,r_0}, L_{b_1}^{p_1),q_1,r_1}\right)_{\theta,r} \cong L_{p,r,a^\#},$$

where $\frac{1}{p} = \frac{1-\theta}{p_0} + \frac{\theta}{p_1}$ and $a^\#(t) = c_1(t)^\theta$.

**Corollary 44.** Let $1 < p_0 < p_1 < \infty$, $0 < r, r_0, r_1, q_1 \leq \infty$, $b_1 \in SV$, $\|s^{-1/r} a(s)\|_{r,(1,\infty)} < \infty$, and $\|s^{-1/r_1} b_1(s)\|_{r_1,(0,1)} < \infty$. Put $c_1(t) = \|s^{-1/r_1} b_1(s)\|_{r_1,(0,t)}$ and $\rho(t) = t^{\frac{1}{p_0}-\frac{1}{p_1}} \frac{1}{c_1(t)}$. Then

$$\left(L_{p_0,r_0}, L_{b_1}^{p_1),q_1,r_1}\right)_{0,r;a} \cong L_{a^\#}^{(p_0,r_0,r)},$$

where $a^\#(t) = a(\rho(t))$.

**Corollary 45.** (Cf. [8, Corollary 7.7].) Let $0 < \theta < 1$, $1 < p_1 < \infty$, $0 < r, r_1, q_1 \leq \infty$, $b_1 \in SV$, and $\|s^{-1/r_1} b_1(s)\|_{r_1,(0,1)} < \infty$. Put $c_1(t) = \|s^{-1/r_1} b_1(s)\|_{r_1,(0,t)}$. Then

$$\left(L_1, L_{b_1}^{p_1),q_1,r_1}\right)_{\theta,r} \cong L_{p,r,a^\#},$$

Where $\frac{1}{p} = 1 - \theta + \frac{\theta}{p_1}$ and $a^\#(t) = c_1(t)^\theta$.

**Corollary 46.** (Cf. [13, Theorem 7.3].) Let $0 < \theta < 1$, $1 < p_0 < p_1 < \infty$, $0 < r, r_0, r_1, q_0 \leq \infty$, $b_0 \in SV$, and $\|s^{-1/r_0} b_0(s)\|_{r_0,(0,1)} < \infty$. Put $c_0(t) = \|s^{-1/r_0} b_0(s)\|_{r_0,(0,t)}$. Then

$$\left(L_{b_0}^{p_0),q_0,r_0}, L_{p_1,r_1}\right)_{\theta,r} \cong L_{p,r,a^\#},$$

where $\frac{1}{p} = \frac{1-\theta}{p_0} + \frac{\theta}{p_1}$ and $a^\#(t) = c_0(t)^{1-\theta}$.

**Corollary 47.** (Cf. [11, Theorem 1.1].) Let $1 < p_0 < p_1 < \infty$, $0 < r, r_0, r_1, q_0 \leq \infty$, $a, b_0 \in SV$, $\|s^{-1/r} a(s)\|_{r,(0,1)} < \infty$, and $\|s^{-1/r_0} b_0(s)\|_{r_0,(0,1)} < \infty$. Put $c_0(t) = \|s^{-1/r_0} b_0(s)\|_{r_0,(0,t)}$ and $\rho(t) = t^{\frac{1}{p_0}-\frac{1}{p_1}} c_0(t)$. Then

$$\left(L_{b_0}^{p_0),q_0,r_0}, L_{p_1,r_1}\right)_{1,r;a} \cong L_{a^\#}^{p_1),r_1,r},$$

where $\frac{1}{p} = \frac{1-\theta}{p_0} + \frac{\theta}{p_1}$ and $a^\#(t) = a(\rho(t))$.

**Corollary 48.** Let $0 < \theta < 1$, $1 < p_0 < \infty$, $0 < r, r_0, q_0 \le \infty$, $b_0 \in SV$, and $\|s^{-1/r_0}b_0(s)\|_{r_0,(0,1)} < \infty$. Put $c_0(t) = \|s^{-1/r_0}b_0(s)\|_{r_0,(0,t)}$. Then
$$\left(L_{b_0}^{p_0),q_0,r_0}, L_\infty\right)_{\theta,r} \cong L_{p,r,a^\#},$$
where $\frac{1}{p} = \frac{1-\theta}{p_0}$ and $a^\#(t) = c_0(t)^{1-\theta}$

### 8.4. Interpolation formulae for couples where both operands are small or grand Lorentz spaces

**Corollary 49.** (Cf. [11, Theorem 3.4].) Let $0 < \theta < 1$, $1 < p_0 < p_1 < \infty$, $0 < r, r_0, r_1, q_0, q_1 \le \infty$, $b_0, b_1 \in SV$, $\|s^{-1/r_0}b_0(s)\|_{r_0,(1,\infty)} < \infty$, and $\|s^{-1/r_1}b_1(s)\|_{r_1,(1,\infty)} < \infty$. Put $c_0(t) = \|s^{-1/r_0}b_0(s)\|_{r_0,(t,\infty)}$ and $c_1(t) = \|s^{-1/r_1}b_1(s)\|_{r_1,(t,\infty)}$. Then
$$\left(L_{b_0}^{(p_0,q_0,r_0)}, L_{b_1}^{(p_1,q_1,r_1)}\right)_{\theta,r} \cong L_{p,r,a^\#},$$
where $\frac{1}{p} = \frac{1-\theta}{p_0} + \frac{\theta}{p_1}$ and $a^\#(t) = c_0(t)^{1-\theta}c_1(t)^\theta$.

*Proof.* Let $\bar{A} = (L_1, L_\infty)$, $\theta_i = 1 - \frac{1}{p_i}$, and $\eta = (1-\theta)\theta_0 + \theta\theta_1$. By Lemma 35 we have
$$L_{b_0}^{(p_0,q_0,r_0)} \cong (L_1, L_\infty)_{\theta_0,r_0,b_0,q_0,1}^{\mathcal{L}} = \bar{A}_{\theta_0,r_0,b_0,q_0,1}^{\mathcal{L}}$$
and
$$L_{b_1}^{(p_1,q_1,r_1)} \cong (L_1, L_\infty)_{\theta_1,r_1,b_0,q_1,1}^{\mathcal{L}} = \bar{A}_{\theta_1,r_1,b_0,q_1,1}^{\mathcal{L}}.$$
Using now Theorem 25 and [6, Corollary 5.3], we get
$$\left(L_{b_0}^{(p_0,q_0,r_0)}, L_{b_1}^{(p_1,q_1,r_1)}\right)_{\theta,r} \cong \left(\bar{A}_{\theta_0,r_0,b_0,q_0,1}^{\mathcal{L}}, \bar{A}_{\theta_1,r_1,b_0,q_1,1}^{\mathcal{L}}\right)_{\theta,r} \cong \bar{A}_{\eta,r,a^\#} \cong L_{p,r,a^\#}.$$
This completes the proof. □

All corollaries below can be proved analogously using theorems from Section 7.

**Corollary 50.** (Cf. [11, Theorem 5.1] and [13, Theorem 6.5].) Let $0 < \theta < 1$, $1 < p_0 < p_1 < \infty$, $0 < r, r_0, r_1, q_0, q_1 \le \infty$, $b_0, b_1 \in SV$, $\|s^{-1/r_0}b_0(s)\|_{r_0,(1,\infty)} < \infty$, and $\|s^{-1/r_1}b_1(s)\|_{r_1,(0,1)} < \infty$. Put $c_0(t) = \|s^{-1/r_0}b_0(s)\|_{r_0,(t,\infty)}$ and $c_1(t) = \|s^{-1/r_1}b_1(s)\|_{r_1,(0,t)}$. Then
$$\left(L_{b_0}^{(p_0,q_0,r_0)}, L_{b_1}^{p_1),q_1,r_1}\right)_{\theta,r} \cong L_{p,r,a^\#},$$
where $\frac{1}{p} = \frac{1-\theta}{p_0} + \frac{\theta}{p_1}$ and $a^\#(t) = c_0(t)^{1-\theta}c_1(t)^\theta$.

**Corollary 51.** (Cf. [11, Theorem 1.2] and [13, Theorem 8.3].) Let $0 < \theta < 1$, $1 < p_0 < p_1 < \infty$, $0 < r, r_0, r_1, q_0, q_1 \le \infty$, $b_0, b_1 \in SV$, $\|s^{-1/r_0}b_0(s)\|_{r_0,(0,1)} < \infty$, and $\|s^{-1/r_1}b_1(s)\|_{r_1,(0,1)} < \infty$. Put $c_0(t) = \|s^{-1/r_0}b_0(s)\|_{r_0,(0,t)}$ and $c_1(t) = \|s^{-1/r_1}b_1(s)\|_{r_1,(0,t)}$. Then
$$\left(L_{b_0}^{p_0),q_0,r_0}, L_{b_1}^{p_1),q_1,r_1}\right)_{\theta,r} \cong L_{p,r,a^\#},$$

where $\frac{1}{p} = \frac{1-\theta}{p_0} + \frac{\theta}{p_1}$ and $a^\#(t) = c_0(t)^{1-\theta} c_1(t)^\theta$.

**Corollary 52.** Let $0 < \theta < 1$, $1 < p_0 < p_1 < \infty$, $0 < r, r_0, r_1, q_0, q_1 \leq \infty$, $b_0, b_1 \in SV$, $\left\| s^{-1/r_0} b_0(s) \right\|_{r_0,(0,1)} < \infty$, and $\left\| s^{-1/r_1} b_1(s) \right\|_{r_1,(1,\infty)} < \infty$. Put $c_0(t) = \left\| s^{-1/r_0} b_0(s) \right\|_{r_0,(0,t)}$ and $c_1(t) = \left\| s^{-1/r_1} b_1(s) \right\|_{r_1,(t,\infty)}$. Then
$$\left( L_{b_0}^{p_0),q_0,r_0}, L_{b_1}^{(p_1,q_1,r_1} \right)_{\theta,r} \cong L_{p,r,a^\#},$$
where $\frac{1}{p} = \frac{1-\theta}{p_0} + \frac{\theta}{p_1}$ and $a^\#(t) = c_0(t)^{1-\theta} c_1(t)^\theta$.